\documentclass[a4paper,11pt]{amsart}
\usepackage[margin=2.7cm]{geometry}
\usepackage{amsfonts,amssymb,amscd,graphicx,stmaryrd,mathtools}
\usepackage[margin=.5in,footnotesize]{caption}
\usepackage[utf8]{inputenc}
\usepackage[T1]{fontenc}
\usepackage[french,english]{babel}
\usepackage[colorlinks]{hyperref}
\usepackage{xcolor,etoolbox}
\usepackage{baskervillef}

\graphicspath{{./fig/}}

\newcommand{\LR}{\mathrm{LR}}

\newcommand{\intint}[2]{\llbracket #1, #2 \rrbracket}

\newtheorem{theorem}{Theorem}

\newtheorem{corollary}[theorem]{Corollary}
\newtheorem{lemma}[theorem]{Lemma}

\newtheorem{proposition}[theorem]{Proposition}

\newtheorem{definition}[theorem]{Definition}
\newtheorem{remark}[theorem]{Remark}

\numberwithin{equation}{section}

\newenvironment{preuve}[1][]
{\begin{proof}[\textbf{Proof#1}]}{\end{proof}}
\newcounter{rk}
\newcommand{\rk}[1]{%
  \stepcounter{rk}%
  \textcolor{red}{\bf \smash{\fbox{XXX \therk}} #1}%
  \csname @latex@warning\endcsname{There are things to do}%
}
\newcounter{assumption} \renewcommand\theassumption{(P\arabic{assumption})}
\newcommand{\newassumption}[1]%
  {\refstepcounter{assumption}\theassumption\label{#1}}

\newcommand{\rkb}[1]{{\color{blue}\bf  #1}}
\newcommand{\rks}[1]{{\color{red} \bf #1}}
\newcommand{\rko}[1]{{\color{orange} \bf #1}}

\newcommand{\bfp}{\mathbf{P}}

\newcommand{\tr}{\mathop{\mathrm{tr}}}

\newcommand{\bt}{\mathcal{T}}
\newcommand{\bv}{\mathcal V}

\newcommand{\Nn}{\mathbb{N}}
\newcommand{\R}{\mathbb{R}}

\newcommand{\Z}{\mathbb{Z}}

\newcommand{\WB}{\mathrm{WB}}

\DeclareMathOperator{\sdiff}{\triangle}

\let\epsilon=\varepsilon%
\let\ln=\log%
\let\leq=\leqslant%
\let\geq=\geqslant%

\author{Vincent Beffara and Damien Gayet}
\title[Percolation without FKG]%
{Percolation without FKG}
\date{\today}

\begin{document}

\begin{abstract}
  We  prove  a  Russo-Seymour-Welsh  theorem  for the antiferromagnetic Ising
  model with small parameter on a regular enough periodic planar triangulation.
  More generally we prove that the box-crossing property extends to a
  neighborhood of the product measure among certain families of discrete
  percolation models that  do not necessarily satisfy the
  Fortuin-Kasteleyn-Ginibre condition  of positive  association.
\end{abstract}

\maketitle

\setcounter{tocdepth}{1}
% \tableofcontents

\section*{Introduction}

Bernoulli site-percolation on a planar periodic lattice is one of the simplest
models of statistical mechanics: for a given parameter $p \in [0,1]$, the
percolation measure $P_p$ is defined by saying that every vertex of the lattice
is \emph{open} with probability $p$ and \emph{closed} with probability $1-p$,
independent of the others. Under minimal assumptions on the lattice, there
exists a critical parameter $p_c$ such that for $p<p_c$ almost surely all the
clusters of open vertices are finite, while for $p>p_c$ there exists with
probability $1$ a unique infinite open cluster.

For a given model, for any rectangle $R\subset \R^2$ and any positive integer
$n$, the \emph{crossing probability} of the scaled rectangle $nR$ is the
probability that there exists an open path in the lattice, contained in $nR$,
and connecting the left and the right side of $nR$; the \emph{box-crossing
property} (BXP for short) states that for every rectangle $R$, the crossing
probability of $nR$ is bounded away from $0$ and $1$ uniformly in $n$. It is a
cornerstone of the study of two-dimensional percolation-like models,
characteristic of critical behavior, and enables for instance to show the
absence of infinite  open clusters at  the critical  point; it  is  the starting
point  of the proofs of existence of conformally invariant scaling limits for
both percolation and the Ising model. It can also be used to derive the
sharpness of phase transition for these models. Russo~\cite{Russo1978} and
Seymour and Welsh~\cite{seymour:rsw} proved in 1976 that Bernoulli percolation
satisfies the BXP at criticality.

All known proofs for the BXP are restricted to models which satisfy a crucial
property, namely the \emph{Fortuin-Kasteleyn-Ginibre  condition} (FKG for
short), meaning that two  increasing events  are positively  correlated: this
allows for instance to construct the pure phases of the Ising model or to build
percolation crossings of long rectangles from more elementary blocks, typically
crossings of squares; it is  an essential  tool in much  of  the  literature in
statistical mechanics. For  Bernoulli percolation, it was first observed by
Harris~\cite{Harris1960}. Köhler-Schindler and Tassion~\cite{kohler2023crossing}
have recently proved that essentially, FKG and natural symmetries are  enough to
have the BXP at an infinite number of scales; in the case of the Ising model,
this applies for all $\beta\geqslant0$.

In the present  paper, we give a novel proof of the box-crossing property which
does not make use of the FKG inequality, at the cost of being restricted to
models that exhibit the domain Markov property and fast enough decorrellation
and which are close enough to being product measures. To  the best of our
knowledge, this is the first proof of such a result in the usual setup of
statistical mechanics, \emph{i.e.} for fixed models  in infinite volume.

A particular case which is of special interest is that of the Ising model. For
every $\beta$ small enough in absolute value, on the face-centered square
lattice there is a unique Gibbs measure $P_\beta$ for the Ising model at inverse
temperature $\beta$, which can be seen as a dependent site-percolation model by
declaring sites with spin $+$ open and sites with spin $-$ closed. Our methods
apply to that setup:
\begin{theorem}\label{thm:main}
  There exists $\beta_0>0$ such that for every $\beta \in (-\beta_0, \beta_0)$,
  the Ising measure $P_\beta$ satisfies the box-crossing property and as a
  consequence, $P_\beta$-almost surely, all the clusters of constant spin are
  finite.
\end{theorem}
The other usual consequences of the box-crossing property, such as th power-law
decay of the connectivity function, the values of the universal exponents, and
the tightness of large interfaces upon scaling limits, should follow by the same
overall arguments as in the case of percolation, though we chose not to include
them in the present paper to focus on the BXP proof itself. It is a natural
direction for future work, with the long-term goal of showing that the
high-temperature Ising model is in the same universality class as Bernoulli
percolation.

The core idea  of the argument is a renormalization scheme at criticality. We
obtain a finite-scale criterion, propagating estimates from  one scale to the
next thanks to a precise control of  correlations rather  than using positive
association; our assumption of fast decorrellation implies the existence of a
critical scale from which good enough estimates propagate to all larger scales.
In turn, the assumed Gibbs property enables estimates at fixed scales to be
obtained by continuity from known results for critical percolation.

\bigskip

\subsection*{Related works}

There is significant literature about Russo-Seymour-Welsh theory and the
box-crossing property in percolation and positively associated models, starting
from the foundational works and Russo~\cite{Russo1978,Russo1981} and Seymour and
Welsh~\cite{seymour:rsw} and culminating in the previously cited article of
Köhler-Schindler and Tassion~\cite{kohler2023crossing}. Much less is known in
the absence of the FKG relation: Beliaev, Muirhead and Wigman~\cite{Beliaev2021}
prove a related estimate for a negatively correllated continuous model but in a
regime where negative association vanishes at large scales;
Vanneuville~\cite{vanneuville} proves a lower bound for the crossing probability
of rectangles at all scales when the model is the sign of a Gaussian field (not
necessarily positively associated), under decay assumptions for the correlation
kernel, though with conditions on the rate of decay that depend on the shape of
the rectangle. In terms of methods, multiscale approaches for critical models
which share similarities to ours were recently used \emph{e.g.}
in~\cite{Severo2020,Severo2023} using a decomposition of the Gaussian free
field.

\subsection*{Open questions}

Our approach raises a number of natural questions. First, given a continuous
one-parameter family of models with appropriate decorrellation, is it the case
that the set of parameters at which the BXP holds is open? In the case of the
Ising model, a very related question would be whether our proof extends beyond a
perturbative regime.  It is often the case that finite-scale arguments allow to
prove similar results; in our case however, the criterion that we use is somehow
not local enough to enable simple use of continuity. It is also natural to
wonder how our assumptions can be relaxed to address a wider range of models,
and in particular the random-cluster model for $q<1$ or continuum models like
level lines of smooth Gaussian fields might be amenable to similar methods. Both
of these directions are beyond what we were able to do.

\subsection*{Acknowledgments} We are deeply grateful to the referee of a
previous (and crucially flawed) version of the paper for their insight. We also
thank Hugo Vanneuville for many useful discussions.

\section{Definitions and setup}

In all this article, we will only consider the face-centered planar square
lattice $\Lambda$ (also known as the \emph{Union-Jack lattice}), embedded in
$\mathbb R^2$ and equipped with the distance induced by $\|\cdot\|_\infty$. Its
set of vertices is the set of points in the plane with either two integer
coordinates or two half-integer coordinates; two vertices are adjacent if and
only if they have distance $\sqrt 2 / 2$, or they have integer coordinates and
have distance $1$. This lattice is invariant under integer translations, as well
as under reflections across coordinate axes, which will be crucial in what
follows. Our construction can be adapted to the case of the usual triangular
lattice, which exhibits a different automorphism group with enough symmetries,
but for sake of readability we will not do this here. Lattices with fewer
symmetries, or which are not triangulations, remain out of reach of our methods;
bond models, such as the usual Bernoulli bond percolation on the square lattice,
exhibit a different version of duality, so it is likely that variations of our
methods would apply to them but we did not attempt to do that here.

We will repeatedly use the notation $\Lambda_R$ for the square box of size $2R$,
namely \[\Lambda_R := (-R,R)^2 \cap \Lambda = \{v \in \Lambda : \|v\|_\infty <
R\},\] and $\mathcal A(r,R)$ for the annulus with inner radius $r$ and outer
radius $R$, defined as \[\mathcal A(r,R) := \{v \in \Lambda : r < \|v\|_\infty <
R\}.\]

\subsection{Probabilistic setup}

In all that follows, if $D$ is a set of vertices, we will denote by $\Omega_D$
the set of configurations on $D$, \emph{i.e.} $\Omega_D = \{\pm1\}^D$. In the
case of the whole lattice we will write $\Omega$ for $\Omega_\Lambda$. When
$\omega$ is an element of $\Omega_D$, $\omega^\ast$ will denote the dual
configuration obtained by applying a global spin-flip, \emph{i.e.} for each $i
\in D$, $\omega^\ast_i = - \omega_i$.

\begin{definition}[Specification]
  Denote by $\mathcal D$ the collection of all pairs $(D, \omega)$ where $D$ is
  a finite set and $\omega \in \Omega_{D^c}$, interpreted as a boundary
  condition. A \emph{specification} is a collection $\mathbb P = (P_{D,
  \omega})_{(D, \omega) \in \mathcal D}$ where for each $(D, \omega) \in
  \mathcal D$, $P_{(D, \omega)}$ is a probability measure on $\Omega_D$,
  satisfying the following compatibility relation: for every $D_1 \subseteq
  D_2$, $\omega \in \Omega_{D_2^c}$ and every event $A$ depending only in the
  configuration within $D_1$, \[ P_{D_2, \omega} [A] = \sum_{\omega' \in
  \Omega_{D_2 \setminus D_1}} P_{D_2, \omega} [{\omega'}] P_{D_1, \omega \cup
  \omega'} [A].\] A probability measure $P$ on $\Omega$ is \emph{compatible}
  with the specification $\mathbb P$ if for every finite set $D$ and for
  $P$-almost every $\omega \in \Omega_{D^c}$, the conditional distribution under
  $P$ of the configuration in $D$ conditionally on $\omega$ is $\mathbb
  P_{D,\omega}$; in particular, for every $D$-measurable event $A$,
  \[\int_{\Omega_{D^c}} P_{D, \omega}[A]\mathrm dP(\omega) = P[A].\]
\end{definition}

All the specifications $\mathbb P = (P_{D, \omega})$ considered in this article
will be assumed to satisfy the following conditions:
\begin{enumerate}
\item[\newassumption{as:unique}] \textbf{Uniqueness:} there is a unique
  probability measure $P$ on $\Omega$ which is compatible with $\mathbb P$;
\item[\newassumption{as:symm}] \textbf{Symmetry:} the specification is invariant
  under lattice automorphisms, in the sense that for every finite $D$, every
  $\omega \in \Omega_{D^c}$ and every automorphism $\theta$ of the lattice,
  $\theta_\ast P_{D, \omega} = P_{\theta(D), \theta_\ast \omega}$.
\item[\newassumption{as:duality}] \textbf{Duality:} the specification is
  invariant under spin-flip, in the sense that for every finite $D$, every
  $\omega \in \Omega_{D^c}$ and every $\sigma \in \Omega_D$, $P_{D, \omega^\ast}
  [\{\sigma^\ast\}] = P_{D, \omega}[\{\sigma\}]$.
\end{enumerate}
We will informally use the term \emph{model} to denote a specification
satisfying conditions \ref{as:unique}--\ref{as:duality}, or the corresponding
probability measure $P$ on $\Omega$. The implementation of the perturbative
argument in Section~\ref{sec:interpolation} will require additional assumptions,
which will be introduced when needed.

As a representative example, the Ising model with inverse temperature $\beta \in
\mathbb R$ corresponds to the specification given by \[ P_{D, \omega}^\beta
[\{\sigma\}] = \frac {\exp (-\beta H_{D}(\sigma;\omega))}{\sum_{\sigma' \in
\{\pm1\}^D} \exp (-\beta H_{D}(\sigma';\omega))}\] for every $\sigma \in
\Omega_D$, where the Hamiltonian $H_{D}$ is defined as \[H_D(\sigma;\omega) = -
\sum_{i,j \in D ; i \sim j} \sigma_i \sigma_j - \sum_{i \in D; j \in D^c; i \sim
j} \sigma_i \omega_j.\] The compatibility relations are easily obtained from the
form of the specification; symmetry and duality are inherited from the
symmetries of the Hamiltonian, and it is classical (see
\emph{e.g.}~\cite{velenik}) that uniqueness holds whenever $|\beta|$ is small
enough.

\subsection{Notation and definitions}

\begin{definition}[Path, Circuit]
  A \emph{path} on $\Lambda$ is a finite sequence $\gamma = (\gamma_0, \ldots,
  \gamma_n)$ of vertices such that for every $i < n$, $\gamma_i$ is adjacent to
  $\gamma_{i+1}$; $n$ is the \emph{length} of $\gamma$, denoted by $|\gamma|$,
  $\gamma_0$ is its \emph{starting point}, and $\gamma_n$ is its
  \emph{endpoint}. A \emph{simple path} is a path with pairwise distinct
  vertices. The \emph{inner part} of $\gamma$ is the set $\mathring\gamma :=
  \{\gamma_k : 0 < k < |\gamma|\}$.

  A \emph{circuit} is a simple path such that its endpoint is adjacent to its
  starting point. The \emph{support} of a circuit is the set of vertices that it
  contains, the \emph{inside} of a circuit is the set of all the vertices that
  its support surrounds, and the \emph{region delimited} by a circuit is the
  union of its support and its interior. A circuit is \emph{direct} if it has
  its interior on its left. Note that the inside of a circuit is not always connected if there is a fjord somewhere on it.
  % \rk{Careful that there can be bubbles, the inside is
  % really the vertices with winding number $1$.} \rk{The region is actually $Q
  % \cup \partial Q$, should be denoted $\bar Q$ or something.} \rk{With this
  % version of the definition, it is not always true that $Q$ is connected.}

  Given two paths $\gamma = (\gamma_0, \ldots, \gamma_n)$ and $\delta =
  (\delta_0, \ldots, \delta_m)$ such that $\gamma_n$ is adjacent to $\delta_0$,
  their \emph{concatenation} is the path $\gamma \oplus \delta := (\gamma_0,
  \ldots, \gamma_n, \delta_0, \ldots, \delta_m)$.
\end{definition}

\begin{definition}[Quad, dual quad]
  A \emph{quad} $Q$ in the lattice $\Lambda$ is a quadruple of paths in which
  the endpoint of each path is adjacent to the starting point of the next, and
  whose concatenation forms a direct simple circuit; equivalently, it can be
  seen as a direct simple circuit partitioned into four non-empty contiguous
  intervals. The four elements of the tuple are the \emph{left (resp.\ bottom,
  right, top) boundaries} of $Q$, and we will denote them by $\partial_\ell Q$,
  $\partial_b Q$, $\partial_r Q$ and $\partial_t Q$ respectively; their
  concatenation is the \emph{boundary} of $Q$, denoted by $\partial Q$.

  The \emph{interior} of the quad $Q$, which will be denoted by $Q^o$, is the
  set of vertices around which $\partial Q$ winds once; the \emph{complement} of
  the quad, denoted by $Q^c$, is the set of vertices that are neither on
  $\partial Q$ nor in $Q^o$. By abuse of notation, the term quad will also be
  used to denote the interior of the circuit $\partial Q$, so we allow ourselves
  to say \emph{e.g.} that a vertex belongs to a quad if it belongs to the $Q^o$.
  Similarly, we will write $v \in \partial_\ell Q$ to mean that the vertex $v$
  belongs to the support of the path $\partial_\ell Q$.

  The \emph{dual} of a quad $Q$ is the quad $Q^\ast$ corresponding to the tuple
  $(\partial_b Q, \partial_r Q, \partial_t Q, \partial_\ell Q)$. It has the same
  support as $Q$ and the left-right and top-bottom boundary pairs swapped (so
  that for instance $\partial_\ell Q^\ast = \partial_b Q$).
\end{definition}

We will use the same notation $\Lambda_R$ for the square defined above, and for
the quad having it as its support and its sides as boundary intervals.

In what follows we will need to control the geometry of various quads, and for
that we will use several metrics which we now define:

\begin{definition}[Distances]
  Let $Q$ be a quad, and let $x$ and $y$ be two vertices belonging to $Q$.
  \begin{enumerate}
  \item The \emph{ambient distance} between $x$ and $y$, denoted by $d(x,y)$, is
    the graph distance in $\Lambda$ between $x$ and $y$, \emph{i.e.} the length
    of the shortest path from $x$ to $y$. It does not depend on $Q$.
  \item The \emph{intrinsic distance} between $x$ and $y$ in $Q^o \cup \partial
    Q$, denoted by $d_Q(x,y)$, is the length of the shortest path $\gamma$ from
    $x$ to $y$ satisfying $\mathring \gamma \subseteq Q^o$ (or $+\infty$ if no
    such path exists).
  \item The \emph{inner distance} between $x$ and $y$ in $Q^o \cup \partial Q$,
    denoted by $\delta_Q(x,y)$, is the smallest \emph{ambient} diameter of a
    path $\gamma$ from $x$ to $y$ satisfying $\mathring \gamma \subseteq Q^o$,
    \emph{i.e.} the smallest $R$ such that there exists a path $\gamma$ from $x$
    to $y$ with $\mathring \gamma \subseteq Q^o$ satisfying \[\forall v, v' \in
    \gamma, \quad d(v,v') \leq R.\]
  \end{enumerate}
\end{definition}

Note that $d_Q$ and $\delta_Q$ do not necessarily satisfy the triangle
inequality: the bound $d_Q(x,z) \leq d_Q(x,y) + d_Q(y,z)$ always holds when $y
\in Q^o$, but \emph{e.g.} in the case when $Q$ is not connected but has two
connected components separated by $y \in \partial Q$, and each containing one of
$x$ and $z$, the left-hand side is infinite while the right-hand side is not.
This will not be an issue for our purpose.

For every $x, y \in Q^o \cup \partial Q$, we have
\[d(x,y) \leq \delta_Q(x,y) \leq d_Q(x,y).\] As usual we extend this notation to
the cases of point-to-set distances and set-to-set distances, defining for
instance for all $x \in Q^o \cup \partial Q$ and $A, B \subseteq Q^o \cup
\partial Q$, \[\delta_Q(x, A) := \inf_{y \in A} \delta_Q(x,y) \quad \text{and}
\quad \delta_Q(A,B) := \inf_{y \in A} \inf_{z \in B} \delta_Q(y,z).\] We will
denote by $D(Q)$ the diameter of $Q$, measured in $\ell^\infty$ ambient
distance, \emph{i.e.} $D(Q)$ is the smallest integer $n$ such that $Q$ fits in a
square box of size $n \times n$.

\begin{definition}[Quad width]
  Let $Q$ be a quad: we define three measures of boundary separations in $Q$,
  \[ W_h(Q) := \delta_Q(\partial_\ell Q, \partial_r Q), \quad W_v(Q) :=
  \delta_Q(\partial_b Q, \partial_t Q) \] and the \emph{width} of $Q$ defined as
  \(W(Q) := \min (W_h(Q), W_v(Q))\).
\end{definition}

\begin{remark}
  Width mixes intrinsic and ambient geometry, so that for every quad $Q$, $W(Q)$
  is bounded below by the smallest ambient distance between vertices on opposite
  sides of $Q$ (which one could term the separation of $Q$), and bounded above
  by the length of the shortest path joining opposite sides of $Q$ (the
  intrinsic width of $Q$); for a ``C''-shaped quad with small opening, the
  separation is much smaller than the width, and if the shape of the ``C'' is
  tortuous the intrinsic diameter can be even larger. It is always true that
  $W(Q) \leq D(Q)$ and $W(Q^\ast) = W(Q)$.
\end{remark}

\begin{definition}[Quad crossings]
  Let $Q$ be a quad.
  \begin{enumerate}
  \item A \emph{(left-to-right, or horizontal) crossing} of $Q$ is a path
    $\gamma$ such that $\gamma_0 \in \partial_\ell Q$, $\gamma_{|\gamma|} \in
    \partial_r Q$, and $\mathring \gamma \subseteq Q$;
  \item We will say that a subset $S$ of $Q$ \emph{crosses $Q$ (horizontally)}
    if there exists a horizontal crossing $\gamma$ of $Q$ such that $\mathring
    \gamma \subseteq S$;
  \item Given a configuration $\omega \in \{\pm1\}^{Q}$, an \emph{open crossing}
    of $Q$ is a crossing $\gamma$ of $Q$ such that for every $x \in \mathring
    \gamma$, $\omega_x = +1$;
  \item We will say that the configuration $\omega$ crosses $Q$ (horizontally)
    if there exists an open horizontal crossing, or equivalently if the set $\{x
    : \omega_x = +1\}$ crosses $Q$ horizontally.
  \end{enumerate}
  We will denote by $\LR(Q)$ the event that there exists an open crossing of
  $Q$.
\end{definition}

% \rk{I changed the above definition because in the constructions below we will
% use explored paths as quad boundaries so the specifications need to extend to
% quad boundaries, but crossings should not use vertices belonging to the boundary
% condition.}

\begin{definition}[Quad ordering]
  Given two quads $Q$ and $Q'$, we say that \emph{$Q$ is shorter than $Q'$}, and
  write $Q \leq Q'$, if every crossing of $Q'$ contains as a subpath a crossing
  of $Q$; equivalently, if the inclusion $\LR(Q') \subseteq \LR(Q)$ holds.
\end{definition}

\begin{definition}[Annulus crossings]
  For every $0 < r < R$, a \emph{(radial) open crossing of the annulus $\mathcal
  A(r, R)$ by a configuration $\omega \in \Omega_{\mathcal A(r,R)}$} is a simple
  path $\gamma$ satisfying $\|\gamma_0\|_\infty=r$,
  $\|\gamma_{|\gamma|}\|_\infty = R$ and $\mathring \gamma \subseteq \{ v \in
  \mathcal A(r,R) : \omega_v=1 \}$.

  Extending the previous notation, if $Q$ is a quad we will write $\mathcal A(r,
  R) \leq Q$ if any open crossing of $Q$ contains as a subpath a radial crossing
  of $\mathcal A(r,R)$; we will denote by $\mathcal Q(r,R)$ the collection of
  all quads for which this condition holds: \[\mathcal Q(r,R) := \{ Q \in
  \mathcal Q : \mathcal A(r, R) \leq Q \}.\] In particular, if $\partial_\ell(Q)
  \subset \Lambda_{r}$ and $\partial_r(Q) \cap \Lambda_R = \varnothing$ then $Q
  \in \mathcal Q(r,R)$.
\end{definition}

% \rk{Examples}

% \begin{definition}[Families of quads]
%   For every $0 < r < R$, let $\mathcal Q(r,R)$ denote the collection of all
%   quads $Q$ forcing a crossing of the annulus $\mathcal A(r,R)$, in the sense
%   that any open crossing of $Q$ contains as a subpath a radial crossing of
%   $\mathcal A(r,R)$.
% \end{definition}

\begin{definition}[Probabilities of crossings]
  Given a quad $Q$ and a specification $\mathbb P$ satisfying the uniqueness and
  symmetry conditions \ref{as:unique} and \ref{as:symm}, we will use three
  related quantities: \(\pi(Q) := P[\LR(Q)]\), \[\quad \pi_-(Q) := \inf_{\omega
  \in \Omega_{Q^c}} P_{Q, \omega}[\LR(Q)] \quad \text{and} \quad \pi_+(Q) :=
  \sup_{\omega \in \Omega_{Q^c}} P_{Q, \omega}[\LR(Q)].\]
  % \rk{I changed the
  % definition to extend the boundary condition to also contain the quad boundary,
  % cf.\ above for the corresponding change to the definition of crossings to
  % match this one.}
  As a shortcut, we will denote by $\pi(w,h)$ the probability
  that there is an open left-to-right crossing of the rectangle $([0, w] \times
  [0, h]) \cap \Lambda$. The following related quantities will be used
  throughout:
  \begin{itemize}
  \item The \emph{oscillation} of a quad, $\delta(Q) := \pi_+(Q) - \pi_-(Q)$;
  \item Their maximum, $\epsilon(R, L) := \max \{ \delta(Q) : W(Q) \geq R \;
    \text{and} \; D(Q) \leq L\}$;
  \end{itemize}
  Last, the central quantity on which our argument relies is defined, for $0 < r
  < R$, as \[\beta(r,R) := \sup_{Q \in \mathcal Q(r, R)} \pi_+(Q).\] When
  considering a family $(\mathbb P_u)$ of models, we will decorate the notation
  accordingly so that \emph{e.g.} $\beta_u(r, R)$ will be defined using $\mathbb
  P_u$ as a probability measure on $\Omega$.
\end{definition}

Notice that $\epsilon$ is non-increasing in its first parameter and
non-decreasing in its second parameter. Notice also that there is no reason for
the existence of a boundary condition realizing $\pi_+$, unless the model enjoys
some kind of domain Markov property: this is an assumption that we will add in
the perturbative argument in the last section of the paper, but which we are not
making here.

\begin{definition}[Good model]
  Let $\mathbb P$ be a model, \emph{i.e.} a specification satisfying conditions
  \ref{as:unique}--\ref{as:duality}, and let $\lambda > 0$: we say that the
  model is \emph{$\lambda$-good up to scale $R$} if it satisfies \[\forall 0 <
  r_1 < r_2 \leq R, \quad \beta(r_1, r_2) \leq \lambda^{-1} \left(\frac {r_1}
  {r_2}\right)^\lambda\] and that it is \emph{$\lambda$-good} (or
  \emph{$\lambda$-good at all scales}) if it satisfies this bound for all $0 <
  r_1 < r_2$. As a shortcut, we will also say that a model is \emph{good} if it is $\lambda$-good for some $\lambda>0$.
\end{definition}

If $\lambda_1 \leq \lambda_2$, any $\lambda_2$-good model is also
$\lambda_1$-good. As an example, Russo-Seymour-Welsh theory shows that Bernoulli
percolation is $\lambda$-good for some $\lambda>0$; assuming conformal
invariance (which is known on the triangular lattice but not on $\Lambda$), it
can be shown that the optimal exponent at large scales is equal to $1/6$
although this does not imply that critical Bernoulli percolation is
$\lambda$-good for any specific $\lambda$ because conformal invariance provides
no information on the multiplicative constant.

\begin{definition}[Measures of model decorrellation]
  If $D$ and $D'$ are sets of vertices such that $D' \subseteq D$, denote by
  $\tilde\alpha(D',D)$ the maximal total variation distance between the
  restriction of $P$ to $D'$, and the restriction to $D'$ of any of the $P_{D,
  \omega}$ for any configuration $\omega$ outside $D$; more formally, \[
  \tilde\alpha (D', D) := \sup_{A \in \sigma(D')} \sup_{\omega \in \Omega_{D^c}}
  |P_{D,\omega}[A] - P[A]|.\] For integers $0<\ell<L$, let
  \[\alpha(\ell,L) := \sup_{\mathrm{diam}(D') \leq L} \sup_{d(D', D^c) \geq
  \ell} \tilde \alpha(D', D).\] We say that a model is
  \emph{$\gamma$-decorrelated}, for $\gamma>0$, if \[\forall L>0, \quad \forall
  0 < \ell < L, \quad \alpha(\ell,L) \leq \gamma^{-1} L^2 \ell^{-2-\gamma}\] and
  that it is \emph{well decorrelated} if it is $\gamma$-decorrelated for some
  $\gamma>0$.
\end{definition}

Bernoulli percolation is clearly well decorrelated, since for it $P_{D,\omega}$
does not actually depend on the boundary condition (\emph{i.e.}, $\tilde\alpha$
is identically zero). As above, a more representative example is the
high-temperature Ising model, for which as we will recall below we in fact get
exponential decorrellation, \emph{i.e.} $\alpha(\ell,L) \leq c^{-1} L^2
e^{-c\ell}$ for some $c>0$.

\section{Topological preliminaries}

We are going to replace FKG with  surgeries of quads and paths for which we will
need a few topological and geometric lemmas. We begin by listing direct
consequences of the definitions above.

\begin{lemma}[Duality]
  Given a quad $Q$ and a configuration $\omega \in \Omega_Q$, either $\omega \in
  \LR(Q)$ or $\omega^\ast \in \LR(Q^\ast)$, but not both. As a consequence, for any pair $(Q, Q')$ of quads, \[Q \leq Q' \iff (Q')^\ast \leq Q^\ast.\]
\end{lemma}

\begin{proof}
  This is classical, but we still include a proof for sake of completeness.

  Assume first that $\omega \in LR(Q)$ and $\omega^\ast \in LR(Q^\ast)$ are both
  true: then there is an open path $\gamma_1$ connecting $\partial_\ell Q$ and
  $\partial_r Q$, and a closed path $\gamma_2$ connecting $\partial_t Q$ and
  $\partial_b Q$. Denote by $a$ and $c$ (resp.\ $b$ and $d$) the endpoints of
  $\gamma_1$ (resp.\ $\gamma_2$). Using $\gamma_1$, $\gamma_2$ and the boundary
  of $Q$, one can form an embedding of the complete graph $K_4$ into the plane
  with vertices at $a$, $b$, $c$ and $d$; since those are on the boundary of
  $Q$, they can also be connected to any fifth point chosen outside $Q$, leading
  to an embedding of $K_5$ into the plane which is impossible.

  Assume now that neither $\omega \in LR(Q)$ nor $\omega^\ast \in LR(Q^\ast)$
  hold. Construct a $3$-coloring of the vertices of $\bar Q$ with the following
  rules:
  \begin{itemize}
  \item Vertices that belong to $\partial_\ell Q$ or are in $Q$, open and
    connected to $\partial_\ell Q$ by an open path have color $1$;
  \item Vertices that belong to $\partial_t Q$ or are in $Q$, closed and
    connected to $\partial_t Q$ by a closed path have color $2$;
  \item All ather vertices have color $3$.
  \end{itemize}
  Note that an open (resp.\ closed) vertex of $Q$ adjacent to a vertex with
  color $1$ (resp.\ $2$) has the same color, so no lattice face with its $3$
  vertices in $Q$ can carry all $3$ colors. Now form a topological triangle with
  edges $\partial_\ell Q$, $\partial_t Q$ and $\partial_b Q \cup \partial_r Q$:
  its $3$ sides are each colored in one of the $3$ colors, so by Sperner's
  lemma, one of the faces of $\bar Q$ has vertices of all $3$ colors. By the
  previous remark, the vertex with color $3$ in that triangle has to lie on
  $\partial_b Q \cup \partial_r Q$, and so either it is on $\partial_b Q$ (in
  which case there is a vertical closed crossing), or it is on $\partial_r Q$
  (in which case there is an open horizontal crossing).

  Finally, assume that $Q \leq Q'$ and let $\omega$ be a configuration in
  $LR(Q^\ast)$: then $\omega^\ast \notin LR(Q)$ by duality, so $\omega^\ast
  \notin LR(Q')$ by assumption, and hence $\omega \in LR((Q')^\ast)$ by duality
  again, thus proving that $LR(Q^\ast) \subseteq LR((Q')^\ast)$.
\end{proof}

\begin{corollary}[Separation iff dual crossing]\label{cor:dualiff}
  Let $Q$ be a quad and let $S$ be a subset of $Q$. Then, $S$ separates
  $\partial_\ell Q$ from $\partial_r Q$ within $Q$, in the usual sense that any
  crossing of $Q$ has to contain a vertex of $S$, if and only if $S$ contains a
  crossing of $Q^\ast$ \emph{i.e.} a path in $Q$ from $\partial_b Q$ to
  $\partial_t Q$.
\end{corollary}

\begin{proof}
  Let $\omega$ be the configuration defined by $\omega(x) = -1$ if $x \in S$ and
  $\omega(x) = 1$ otherwise. Then $S$ separates $\partial_\ell Q$ from
  $\partial_r Q$ if and only if any crossing contains a closed vertex for
  $\omega$, which is equivalent to the fact that $\omega$ does not belong to
  $LR(Q)$. By the previous lemma, this is equivalent to saying that $\omega^\ast
  \in LR(Q^\ast)$ which is exactly the event that $\omega^\ast$ crosses
  $Q^\ast$, meaning that there is a crossing of $Q^\ast$ that is contained in
  $\{x : \omega^\ast (x) = 1\}$ which is exactly $S$.
\end{proof}

\begin{lemma}[Basic bounds]
  \label{lem:obviousbounds}
  For every quad $Q$, \[ \pi_-(Q) \leq \pi(Q) \leq \pi_+(Q). \] More generally,
  if $D$ is a finite vertex set and $Q \subseteq D$, then for every $\omega \in
  \Omega_{D^c}$, \[\pi_-(Q) \leq P_{D,\omega}[\LR(Q)] \leq \pi_+(Q).\]
\end{lemma}

\begin{proof}
  The first point follows directly from the definitions: denoting by $E_{Q^c}$
  the expectation over $\omega \in \Omega_{Q^c}$ with the induced law obtained
  by restricting $P$, one has \[ \pi(Q) = E_{Q^c}[P_{Q,\omega}[\LR(Q)]] \geq
  \inf_{\omega \in \Omega_{Q^c}} P_{Q,\omega}[\LR(Q)] = \pi_-(Q) \] and
  similarly for the uper bound. For the second one, fix $Q \subseteq D$ and
  $\omega \in \Omega_{D^c}$: one has \[ P_{D,\omega}[\LR(Q)] = \sum_{\eta \in
  \Omega_{D \setminus Q}} P_{D,\omega}[\{\eta\}] P_{Q, \omega \cup \eta}[\LR(Q)]
  \geq \sum_{\eta \in \Omega_{D \setminus Q}} P_{D,\omega}[\{\eta\}] \pi_-(Q) =
  \pi_-(Q) \] where as usual $\omega \cup \eta$ denotes the configuration on
  $Q^c$ which agrees with $\omega$ on $D^c$ and with $\eta$ on $D \setminus Q$;
  and again the upper bound is exactly similar.
\end{proof}

\begin{lemma}[Estimates from duality]
  For every quad $Q \subseteq D$ and every boundary condition $\omega$, \[
  P_{D,\omega} [\LR(Q)] + P_{D,\omega^\ast} [\LR(Q^\ast)] = 1 \quad \text{and}
  \quad \pi(Q) + \pi(Q^\ast) = 1. \] In particular, $\pi(\Lambda_R) = 1/2$.
\end{lemma}

\begin{lemma}[Ordering]
  \label{lem:ordering}
  For any two quads $Q \leq Q'$, the following hold:
  \begin{enumerate}
  \item $W_h(Q) \leq W_h(Q')$, $\pi(Q) \geq \pi(Q')$ and $\pi_-(Q) + \delta(Q)
    \geq \pi_+(Q') - \delta(Q')$.
  \item If in addition $Q \subseteq Q'$, then $\pi_+(Q) \geq \pi_+(Q')$.
  \end{enumerate}
  Moreover, for every $r' \leq r \leq R \leq R'$, $\beta(r,R) \geq
  \beta(r',R')$.
\end{lemma}

\begin{proof}
  (1) Let $\gamma$ be a horizontal crossing of $Q'$ minimizing ambient diameter,
  hence $\mathrm{diam} (\gamma) = W_h(Q')$. By assumption, there exists a
  sub-path $\gamma' \subseteq \gamma$ crossing $Q$ horizontally, so \[ W_h(Q)
  \leq \mathrm{diam} (\gamma') \leq \mathrm{diam} (\gamma) = W_h(Q').\] The
  other inequalities follow directly from the definition of relation $Q \leq
  Q'$.

  (2) Pick $\eta>0$ and let $\omega$ be a boundary condition outside $Q'$
  realizing $P_{Q',\omega} \geqslant \pi_+(Q')-\eta$: then, \[\pi_+(Q) \geq
  P_{Q',\omega}[\LR(Q)] \geq P_{Q',\omega}[\LR(Q')] \geqslant \pi_+(Q') -
  \eta,\] and since this holds for all $\eta>0$ we get $\pi_+(Q) \geq
  \pi_+(Q')$.

  To conlude: for every quad $Q \in \mathcal Q(r',R')$ it clearly holds that $Q
  \in \mathcal Q(r,R)$. From the previous point, we hence get $\pi_+(Q) \leq
  \beta(r,R)$. Taking the supremum over $Q$ concludes the proof.
\end{proof}

Note that in point~(2) the inclusion $Q \subseteq Q'$ is necessary for the bound
to hold in all generality.

\begin{lemma}[Basic properties of distances in quads] Let $Q$ be a quad: for
  every vertex $x \in Q$, \[d(x, \partial Q) = d_Q(x, \partial Q) = \delta_Q(x,
  \partial Q).\]
\end{lemma}

\begin{proof}
  (1) Clearly $d(x, \partial Q) \leq \delta_Q(x, \partial Q) \leq d_Q(x,
  \partial_Q)$ by comparison between the distances, so we only need to show that
  $d_Q(x, \partial Q) \leq d(x, \partial Q)$; but the closest point to $x$ along
  $\partial Q$ in ambient distance, say $y$, is connected to $x$ by a path
  $\gamma$ of length $d(x, \partial Q)$ along which by minimality no other point
  of $\partial Q$ can lie, so $\gamma \subseteq Q$ and hence $d_Q(x, \partial Q)
  \leq d_Q(x,y) \leq d(x, \partial Q)$.
\end{proof}

\begin{lemma}[Quad shaving]
  \label{lem:shaving}
  Let $Q$ be a quad, and let $0 < \ell < W(Q)/2$. Then there exists a quad $Q'
  \subseteq Q$ satisfying the following two conditions:
  \begin{enumerate}
    \item $\forall x \in Q', d(x, \partial Q) \geq \ell$;
    \item $\forall x \in \partial_\ast Q', \delta_Q(x, \partial_\ast Q) = \ell$
    (where $\ast$ is any of $\ell$, $r$, $t$ or $b$).
  \end{enumerate}
\end{lemma}

\begin{proof}
  We proceed in four steps, treating each side of $Q$ in order. First, let \[K_1
  := \{ x \in Q : \delta_Q (x, \partial_l Q) = \ell \}.\] Any path crossing $Q$
  from left to right starts at intrinsic distance $0$ from $\partial_l Q$, and
  ends at intrinsic distance to it at least equal to $W_h(Q) > \ell$, so it has
  to pass through a vertex lying exactly at intrinsic distance $\ell$ from
  $\partial_l Q$, hence it has to intersect $K_1$. This shows that $K_1$
  separates $\partial_l Q$ from $\partial_r Q$, and hence, by
  Corollary~\ref{cor:dualiff}, contains a simple path $\Gamma_1$ crossing $Q$
  vertically. Let $Q_1$ be the quad having $\Gamma_1$ as its left boundary,
  $\partial_r Q$ as its right boundary, and the portions of $\partial_t Q$ and
  $\partial_b Q$ between $\Gamma_1$ and $\partial_r Q$ as top and bottom
  boundaries, respectively. Every vertex in $Q_1$ is separated from $\partial_l
  Q$ by $\Gamma_1$ and therefore lies at $\delta_Q$-distance at least $\ell$
  from $\partial_l Q$. Besides, let $\gamma$ be a horizontal crossing of $Q_1$
  with minimal ambient diameter $W_h(Q_1)$: $\gamma$ can be extended by a path
  of length $\ell$, and therefore diameter at most $\ell$, to form a horizontal
  crossing of $Q$, thus proving that $W_h(Q) \leq W_h(Q_1) + \ell$ or in other
  words $W_h(Q_1) \geq W_h(Q) - \ell$.

  Now let $K_2 := \{ x \in Q_1 : \delta_Q (x, \partial_t Q) = \ell\}$ (note that
  the distance is still measured relative to the original quad $Q$ and not
  within $Q_1$). As above, $K_2$ separates $\partial_t Q_1$ from $\partial_b
  Q_1$ within $Q_1$ so it contains a simple path $\Gamma_2$ crossing $Q_1$
  horizontally. Let $Q_2$ be the sub-quad of $Q_1$ having $\Gamma_2$ as its top
  boundary. Every vertex in $Q_2$ is at $\delta_Q$-distance at least $\ell$ from
  $\partial_t Q$ --- and also from $\partial_l Q$ because $Q_2 \subseteq Q_1$.
  As above, $W_v(Q_2) \geq W_v(Q) - \ell$, and by inclusion, $W_h(Q_2) \geq
  W_h(Q_1) \geq W_h(Q) - \ell$.

  Then, let $K_3 := \{ x \in Q_2 : \delta_Q (x, \partial_r Q) = \ell \}$. We
  still have $\ell < W_h(Q_2)$ so $K_3$ separates $\partial_l Q_2$ from
  $\partial_r Q_2$ within $Q_2$, and it contains a vertical simple crossing
  $\Gamma_3$ of $Q_2$. Let $Q_3$ denote the sub-quad of $Q_2$ having $\Gamma_3$
  as its right boundary. Every vertex in $Q_3$ is at $\delta_Q$-distance at
  least $\ell$ from all of $\partial_l Q$, $\partial_t Q$ and $\partial_r Q$.
  Moreover, $W_v(Q_3) \geq W_v(Q) - \ell$ and $W_h(Q_3) \geq W_h(Q) - 2\ell$.

  Last, let $K_4 := \{ x \in Q_3 : \delta_Q (x, \partial_b Q) = \ell \}$. As
  before it separates $\partial_t Q_3$ from $\partial_b Q_3$, so it contains a
  simple horizontal simple crossing $\Gamma_4$ of $Q_3$: let $Q'$ denote the
  sub-quad of $Q_3$ lying above $\Gamma_4$. It remains to show that $Q'$
  satisfies both needed contitions. Condition (2) follows directly from the
  definition of the $K_i$; and from the construction, every vertex in $Q'$ is at
  $\delta_Q$-distance at least $\ell$ from every side of $Q$, hence from
  $\partial Q$ itself. But since we know that for every $v \in Q$, $d(v,
  \partial Q) = \delta_Q(v, \partial Q)$, this ensures that condition (1) is
  verified as well.
\end{proof}

\section{RSW steps in terms of \texorpdfstring{$\epsilon$}{ε}}
\label{sec:rsw}

In this section, we follow the usual proof of Russo-Seymour-Welsh estimates from
percolation and replace the Harris-FKG inequality everywhere with random quads
and oscillation estimates. The structure of the proof is therefore quite
standard, and our error measure $\varepsilon(R,L)$ is taylored for use here. In
all that follows, for $\lambda_1, \lambda_2 > 1$ we will denote by $\mathcal
L(R, \lambda_1, \lambda_2)$ the L-shaped quad with region
\[\llbracket 0, \lambda_1R\rrbracket \times \llbracket 0,R \rrbracket \; \cup \;
\llbracket (\lambda_1-1)R, \lambda_1 R \rrbracket \times \llbracket 0,
\lambda_2R \rrbracket,\] left boundary $\{0\} \times \llbracket 0,R \rrbracket$
and right boundary $\llbracket (\lambda_1-1)R, \lambda_1R \rrbracket \times
\{\lambda_2 R\}$, and by $\mathcal A(R)$ the annulus $\mathcal A(R, 2R)$ (which
can be covered with two isometric copies of $\mathcal L(R, 4, 4)$).

\begin{lemma}[Rectangle to L]
  \label{lem:rect_to_L}
  For every $\lambda_1, \lambda_2 \geq 3/2$ and every $R$,
  \[\pi(\mathcal L(R, \lambda_1, \lambda_2)) \geq \pi(\lambda_1R, R)
  \pi(\lambda_2R, R) - \varepsilon(R/2, \max(\lambda_1,\lambda_2)R).\]
\end{lemma}

\begin{proof}
  Denote $R_1 := \llbracket 0, \lambda_1R\rrbracket \times \llbracket 0,R
  \rrbracket$ and $R_2 := \llbracket (\lambda_1-1)R, \lambda_1 R \rrbracket
  \times \llbracket 0, \lambda_2R \rrbracket$ the two rectangles forming
  $\mathcal L(R, \lambda_1, \lambda_2)$. On the event $\LR(R_1)$, denote by
  $\Gamma$ the lowest horizontal crossing of $R_1$ and by $Q_\Gamma$ the
  connected component of $(R_1 \cup R_2) \setminus \Gamma$ having
  $\partial_r(\mathcal L(R, \lambda_1, \lambda_2))$ on its boundary (in grey on
  Figure~\ref{fig:rectangle_to_L}); equip it with a quad structure by letting
  $\partial_\ell(Q_\Gamma) = \Gamma$ and $\partial_r(Q_\Gamma) =
  \partial_r(\mathcal L(R, \lambda_1, \lambda_2))$.

  \begin{figure}[ht]
    \includegraphics{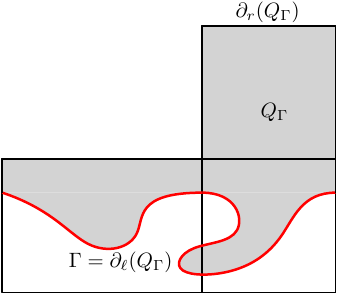}
    \caption{The random quad $Q_\Gamma$ used in the proof of
    Lemma~\ref{lem:rect_to_L}.}
    \label{fig:rectangle_to_L}
  \end{figure}

  By construction $Q_\Gamma \leq R_2$, $D(Q_\Gamma) \leq \max(\lambda_1,
  \lambda_2) R$ and $W(Q_\Gamma) \geq R/2$, hence if $\omega$ is a boundary
  condition outside $L$:
  \begin{align*}
    P_\omega(\LR(\mathcal L(R, \lambda_1, \lambda_2))) &\geq P_\omega(\LR(\mathcal L(R, \lambda_1, \lambda_2)) \cap \LR(R_1)) \\
    &= \sum_{\gamma} P_\omega (\LR(\mathcal L(R, \lambda_1, \lambda_2)) \cap \LR(R_1) ; \Gamma = \gamma) \\
    &\geq \sum_\gamma P_\omega [\LR(R_1) ; \Gamma = \gamma ; \LR(Q_\gamma)] \\
    &\geq \sum_{\gamma} P_\omega(\LR(R_1) ; \Gamma = \gamma) \pi_-(Q_\gamma) \\
    &\geq \sum_\gamma P_\omega(\LR(R_1) ; \Gamma = \gamma) (\pi(Q_\gamma) - \delta(Q_\gamma)) \\
    &\geq \sum_\gamma P_\omega(\LR(R_1) ; \Gamma = \gamma) (\pi(R_2) - \epsilon(R/2, \max(\lambda_1,\lambda_2)R)) \\
    &= P_\omega [\LR(R_1)] (\pi(R_2) - \epsilon(R/2, \max(\lambda_1,\lambda_2)R)) \\
    &\geq (\pi(R_1) - \epsilon(R, \lambda_1 R)) (\pi(R_2) - \epsilon(R/2, \max(\lambda_1,\lambda_2)R)).
  \end{align*}
  The conclusion follows using the monotonicity of $\epsilon$ in its parameters
  and the fact that $\pi(R_1)$ and $\pi(R_2)$ are both bounded above by $1/2$.
\end{proof}

\begin{lemma}[Long rectangle to longer rectangle]
  \label{lem:longer}
  For every $m \geq 3/2$ and $R>0$,
  \[ \pi ((2m-1)R , R) \geq
  \frac {\pi(mR,R)^2}{2} - \epsilon(R/2, mR).\]
\end{lemma}

\begin{proof}
  The proof is similar to that of Lemma~\ref{lem:rect_to_L}. Fix $m \geq 3/2$
  and $R > 0$; let $Q_1$ be the rectangle $\llbracket 0, mR\rrbracket \times
  \llbracket 0, R\rrbracket$ and $Q_2$ be the square $\llbracket (m-1)R,
  mR\rrbracket \times \llbracket 0, R\rrbracket$, both seen as quads with their
  horizontal orientation. Let $A$ be the event that there exists a connected
  cluster inside $Q_1$ which crosses $Q_2$ vertically and touches $\partial_\ell
  Q_1$; denote by $\Gamma$ the rightmost vertical crossing of $Q_2$ (which
  exists on $A$). For every simple vertical crossing $\gamma$ of $Q_2$, denote
  by $Q_\gamma$ the quad formed by the part of $Q_1$ lying to the left of
  $\gamma$, with $\partial_r R_\gamma = \gamma$. Note that $Q_\gamma \leq Q_1$
  and $W(Q_\gamma) \geq R/2$. Then,
  \begin{align*}
    P[A] = \sum_\gamma P[A ; \Gamma = \gamma] &= \sum_\gamma P[\LR(Q_2^\ast) ; \Gamma = \gamma] P[A | \LR(Q_2^\ast) ; \Gamma = \gamma] \\
    &\geq \sum_\gamma P[\LR(Q_2^\ast) ; \Gamma = \gamma] \pi_-(Q_\gamma) \\
    &\geq \sum_\gamma P[\LR(Q_2^\ast) ; \Gamma = \gamma] (\pi(Q_\gamma) - \delta(R_\gamma)) \\
    &\geq \sum_\gamma P[\LR(Q_2^\ast) ; \Gamma = \gamma] (\pi(Q_1) - \varepsilon(R/2, mR)) \\
    &= (\pi(mR,R) - \varepsilon(R/2, mR)) \sum_\gamma P[\LR(Q_2^\ast) ; \Gamma = \gamma] \\
    &= (\pi(mR,R) - \varepsilon(R/2, mR)) P[\LR(Q_2^\ast)] \\
    &= \frac {\pi(mR,R) - \varepsilon(R/2, mR)} {2}.
  \end{align*}

  Now, on the event $A$ one can choose the leftmost vertical crossing $\Gamma'$
  of $Q_2$ that is part of the (unique) cluster realizing $A$; it is measurable
  in the state of the confirugation along $\Gamma$ and to its left, and delimits
  a quad to its right which is shorter than $Q_3 := \llbracket (m-1)R, (2m-1)R
  \rrbracket \times \llbracket 0,R \rrbracket$, so applying the same constuction
  as above to the right of $\Gamma'$ we get
  \begin{align*}
    \pi((2m-1)R, R) &\geq P[A] (\pi(mR, R) - \varepsilon(R/2, mR)) \\
    &\geq \frac
    {(\pi(mR,R) - \varepsilon(R/2, mR))^2} {2} \geq \frac {\pi(mR,R)^2}{2} -
    \epsilon(R/2, mR)
  \end{align*}
  concluding the proof of the lemma.
\end{proof}

\begin{remark}
  In fact, the proof above shows a slightly stronger result, which is the same lower bound on
  \[P \big[\LR(\intint{0}{(2m-1)R} \times \intint{0}{R}) \cap
  \LR((\intint{(m-1)R}{mR} \times \intint{0}{R})^\ast)\big] \] \emph{i.e.} on
  the probability that the ``long rectangle'' is crossed horizontally while at
  the same time the middle square is crossed vertically. Indeed, the cluster
  built as a witness of the event always contains a vertical crossing of $Q_2$
  as part of the scaffolding.
\end{remark}

\begin{lemma}[Half quad]
  \label{lem:half}
  Let $w>0$ and let $Q$ be a quad which is symmetric across the horizontal axis
  and satisfies $W(Q) \geq w$; let $Q/2$ denote the quad with the same support
  but right boundary set to $\partial_r(Q/2) := \partial_r(Q) \cap (\mathbb R
  \times \mathbb R_-)$, and assume that $\mathrm d(\partial_r(Q/2),
  \partial_t(Q)) \geq w/2$. Then $\pi(Q/2) \geq \pi(Q)/2$ and
  \[ \pi_-(Q/2) \geq \frac {\pi_-(Q)} 2 - \epsilon(w/2, D(Q)).\]
\end{lemma}

\begin{proof}
  The bound for $\pi$ it is just a union bound using the symmetry of the model.
  For $\pi_-$, notice that under the geometric assumptions on $Q$, we have
  $W(Q/2) \geq w/2$ so we can write $\pi_-(Q/2) \geq \pi(Q/2) -
  \epsilon(w/2, D(Q))$.
\end{proof}

\begin{remark}
  Lemma~\ref{lem:half} applies in particular to symmetrical quads with $W(Q)
  \geq w$ satisfying the inclusion $\partial_t(Q) \subseteq \mathbb R \times
  \{w/2\}$, for which the condition $\mathrm d(\partial_r(Q/2), \partial_t(Q))
  \geq w/2$ is obvious.
\end{remark}

\begin{lemma}[Square to Rectangle]
  For every $R>0$, \[\pi(3R/2, R) \geq \frac 1 {32} - \epsilon(R/2, R).\]
\end{lemma}

\begin{proof}
  Let $Q_1$ be the square quad $\llbracket R/2, R \rrbracket \times \llbracket
  0, R/2  \rrbracket$, and denote by $\Gamma$ the rightmost vertical crossing of
  it (if it exists, which happens with probability $1/2$). Let $S$ be the
  reflection across the horizontal line $\mathbb R \times \{R/2\}$; if $\gamma$
  is a vertical crossing of $Q_1$, let $Q_\Gamma$ be the quad with right
  boundary $\gamma \cup S(\gamma)$, left boundary $\{0\} \times \llbracket 0, R
  \rrbracket$ and horizontal boundaries contined in those of $\llbracket 0,
  R\rrbracket^2$. From the previous lemma, for every such $\gamma$, the
  probability that there is an open path from $\partial_\ell Q_\gamma$ to
  $\gamma$ contained within $Q_\gamma$ is at least equal to $1/4 -
  \varepsilon(R/2,R)$. Now let $A$ be the event that there exists an open
  crossing $\Gamma$ of $Q_1$ connected within $Q_\Gamma$ to $\{0\} \times
  \llbracket 0, R \rrbracket$ : summing over all possible vertical crossings of
  $Q_1$ like in the proof of Lemma~\ref{lem:longer}, and using the previous
  estimate, we get
  \[P[A] \geq 1/2 (1/4 - \varepsilon(R/2,R)) \geq 1/8 - \varepsilon(R/2,R).\]

  Proceeding again as in the proof of Lemma~\ref{lem:longer} and using
  Lemma~\ref{lem:half} a second time, we obtain that conditionally on $A$, the
  probability that the \emph{rightmost} $\Gamma$ realizing it is connected
  within $\llbracket 0, 3R/2 \rrbracket \times \llbracket 0,R \rrbracket$ to its
  right boundary is at least equal to $1/4 - \varepsilon(R/2, R)$. If this
  occurs, the quad $\llbracket 0, 3R/2 \rrbracket \times \llbracket 0,R
  \rrbracket$ is crossed horizontally, so we get the lower bound \[\pi(3R/2, R)
  \geqslant \left( \frac18 - \varepsilon(R/2, R) \right) \left( \frac14 -
  \varepsilon(R/2, R)\right) \geq \frac 1 {32} - \varepsilon(R/2,R)\] concluding
  the proof of the lemma.
\end{proof}

In what follows, let $\mathcal L'(R)$ be the event (``decorated L-crossing'')
that all the following happen:
\begin{enumerate}
\item The L-shaped quad $\mathcal L(R,4,4)$ is crossed;
\item The square $\intint{R}{2R} \times \intint{0}{R}$ is crossed vertically;
\item The square $\intint{3R}{4R} \times \intint{2R}{3R}$ is crossed
  horizontaly.
\end{enumerate}

\begin{lemma}[Decorated L]
  For every $R\geq 1$,
  \[P[\mathcal L'(R)] \geq \frac {\pi(2R,R)^2 \pi(\mathcal L(R,3,3))} {4} -
  \varepsilon(R/2,4R).\]
\end{lemma}

\begin{proof}
  Let $A$ be the event that there is a vertical crossing of $\intint{R}{2R}
  \times \intint{0}{R}$ which is connected to $\{0\} \times \intint{0}{R}$
  within $\intint{0}{2R} \times \intint{0}{R}$, and $B$ be the event that there
  is a horizontal crossing of $\intint{3R}{4R} \times \intint{2R}{3R}$ which is
  connected to $\intint{3R}{4R} \times \{4R\}$ within $\intint{3R}{4R} \times
  \intint{2R}{4R}$. $A$ and $B$ have the same probability which is bounded
  below, from the argument in the proof of Lemma~\ref{lem:longer}, by
  $\pi(2R,R)/2 - \varepsilon(R/2,2R)$. Moreover, by the same proof,
  \begin{align*}
  P[B|A] & \geq \pi_-(R,R) (\pi(2R,R) - \varepsilon(R/2,2R)) \\
  & \geq (1/2 - \varepsilon(R,R))(\pi(2R,R)
  - \varepsilon(R/2,2R)) \geq \frac {\pi(2R,R)}{2} - \varepsilon(R/2,2R),
  \end{align*}
  so $P[A \cap B] \geq \pi(2R,R)^2/4 - \varepsilon(R/2,2R)$.

  Now on the event $A \cap B$, the leftmost witness of $A$ and the highest
  witness of $B$ delimit a random quad which is shorter than an isometric image
  of $\mathcal L(R,3,3)$, and has diameter less than $4R$ and width at least
  $R$; if that quad is then crossed, then $\mathcal L'(R)$ occurs, so we obtain,
  again summing over possible realizations,
  \begin{align*}
  P[\mathcal L'(R)] & \geq P[A\cap B]
  (\pi(\mathcal L(R,3,3)) - \varepsilon(R,4R)) \\
  & \geq \left(\frac {\pi(2R,R)^2} {4} -
  \varepsilon(R/2,2R) \right) (\pi(\mathcal L(R,3,3)) - \varepsilon(R,4R))
  \end{align*}
  which, by monotonicity of $\epsilon$ in its arguments, concludes the proof.
\end{proof}

\begin{lemma}[L to annulus]
  \label{lem:L-to-annulus}
  \[\pi_-[\mathcal A(R)] \geq P[\mathcal L'(R)] \pi(\mathcal L(R, 4, 4)) -
  \epsilon(R/2,4R).\]
\end{lemma}

\begin{proof}
  The argument is conceptualy very similar to the proof of
  Lemma~\ref{lem:rect_to_L}, with the notable difference that the gluing has to
  be done simultaneously at both ends of a random quad which requires some
  caution. Let $Q_{NW} := \intint{-2R}{-R} \times \intint{-2R}{2R} \;\cup\;
  \intint{-2R}{2R} \times \intint{R}{2R}$ denote the quad forming the
  ``north-west corner'' of $\mathcal A(R)$, which is isometric to $\mathcal
  L(R,4,4)$, and similarly let $Q_{SE} := \intint{-2R}{2R} \times
  \intint{-2R}{-R} \;\cup\; \intint{R}{2R} \times \intint{-2R}{2R}$ denote the
  south-east corner. We will denote by $E$ the event that $Q_{NW}$ is crossed together with two transversal squares, like in $\mathcal L'(R)$.

  \begin{figure}[ht]
    \centering
    \includegraphics{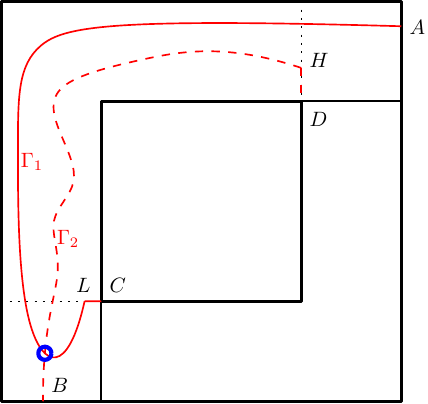}
    \caption{The contradiction in the topological argument in
      Lemma~\ref{lem:L-to-annulus}.}
    \label{fig:annulus2}
  \end{figure}

  Let $\Gamma$ be the outermost simple crossing of $Q_{NW}$, if it exists; let
  $A$ be the point of $\Gamma$ on $\partial_r (Q_{NW})$, let $B$ be its point on
  $\partial_\ell (Q_{NW})$, and let $H$ (resp.\ $L$) be its lowest (resp.\
  rightmost) intersection with the segment $\{R\} \times \intint{R}{2R}$ (resp.\
  $\intint{-2R}{-R} \times \{-R\}$). We first claim that $H$ is visited before
  $L$ on $\Gamma$ going from $A$ to $B$. This can be shown by contradiction:
  assume that instead $H$ is visited after $L$, and let then $\Gamma_1$ denote
  the portion of $\Gamma$ from $A$ to $L$, with the segment $[LC]$ appended, and
  $\Gamma_2$ is the portion of $\Gamma$ from $H$ to $B$, with the segment $[HD]$
  appended; see Figure~\ref{fig:annulus2}. These are two disjoint curves joining
  respectively $A$ to $C$ and $B$ to $D$ within $Q_{NW}$, but because of the
  cyclic order of the $4$ points along the boundary, they have to intersect,
  leading to a contradiction.

  \begin{figure}[ht]
    \centering
    \includegraphics{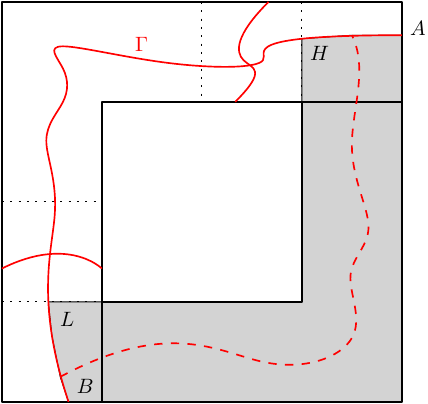}
    \caption{The setup of the topological argument in
      Lemma~\ref{lem:L-to-annulus}.}
    \label{fig:annulus}
  \end{figure}

  Now, assume $E$ is satisfied and as in Figure~\ref{fig:annulus}, define $Q$ to
  be the quad having as left boundary the portion of $\Gamma$ from $L$ to $B$,
  as right boundary the posrion of $\Gamma$ from $A$ to $H$ (those two being
  disjoint by assumption), and as top and bottom boundaries the portions of the
  top and bottom boundaries of $Q_{SE}$ going respectively from $L$ to $H$ and
  from $B$ to $A$. (This quad is shaded in gray in the figure.) $Q$ is shorter
  than $Q_{SE}$ (though not necessarily contained in it), and it satisfies $D(Q)
  \le 4R$ and $W(Q) \geq R$, so conditionally on the configuration outside it,
  it is crossed with conditional probability at least $\pi(\mathcal L(R,4,4)) -
  \varepsilon(R,4R)$; and besides, any crossing of it can be completed with a
  portion of $\Gamma$ to form a circuit surrounding the annulus.

  Summing as before over all possible realizations of the first case, we get
  \[\pi_-(\mathcal A(R)) \geq (P[E] - \varepsilon(R/2,4R)) (\pi(\mathcal
  L(R,4,4)) - \varepsilon(R,4R))\]
  from which the claim again follows by monotonicity of $\epsilon$.
\end{proof}

\begin{theorem}[Annulus estimate in terms of $\varepsilon$]
  \label{thm:annulus}
  For every $R>0$, the probability that the annulus $\mathcal A(R, 2R)$ contains
  an open circuit satisfies
  \[\pi(\mathcal A(R,2R)) \geq \pi_-(\mathcal A(R,2R)) \geq 2^{-140} -
  2\epsilon(R/2, 4R).\]
\end{theorem}

\begin{proof}
  This is a direct consequence of all the previous lemmas. Set $\tilde\epsilon
  := \epsilon(R/2,4R)$: in order,
  \begin{align*}
    \pi(3R/2,R) &\geq 2^{-5} - \varepsilon(R/2, R) \geq 2^{-5} - \tilde\epsilon\\
    \pi(2R,R) &\geq \frac{\pi(3R/2,R)^2}{2} - \varepsilon(R/2,3R/2) \geq \frac {(2^{-5}
      - \tilde\epsilon)^2}{2} - \tilde\epsilon \geq 2^{-11} - 2\tilde\epsilon \\
    \pi(3R,R) &\geq \frac {\pi(2R,R)^2} {2} - \epsilon(R/2,2R) \geq \frac {(2^{-11}
      - 2\tilde\epsilon)^2} {2} - \tilde\epsilon \geq 2^{-23} - 2\tilde\epsilon \\
    \pi(5R,R) &\geq \frac {\pi(3R,R)^2} {2} - \varepsilon(R/2,3R) \geq \frac {(2^{-23}
      - 2\tilde\epsilon)^2} {2} - \tilde\epsilon \geq 2^{-47} - 2\tilde\epsilon \\
    \pi(\mathcal L(R,3,3)) &\geq \pi(3R,R)^2 - \epsilon(R/2,3R) \geq (2^{-23} -
      2\tilde\epsilon)^2 - \tilde\epsilon \geq 2^{-46} - 2\tilde\epsilon \\
    P[\mathcal L'(R)] &\geq \frac {\pi(2R,R)^2 \pi(\mathcal L(R,3,3))} {4} -
      \tilde\epsilon \geq \frac {(2^{-11} - 2\tilde\epsilon)^2 (2^{-46} -
      2\tilde\epsilon)} {4} - \tilde\epsilon \\
    &\geq 2^{-70} - 2 \tilde\epsilon \\
    \pi(\mathcal L(R,4,4)) &\geq P[\mathcal L'(R)] \geq 2^{-70} - 2 \tilde\epsilon \\
    \pi_-(\mathcal A(R,2R)) &\geq P[\mathcal L'(R)]\pi(\mathcal L(R,4,4)) -
      \tilde\epsilon \geq (2^{-70} - 2 \tilde\epsilon)^2 - \tilde\epsilon
  \end{align*}
  concluding the proof.
\end{proof}

\section{Regularity of \texorpdfstring{$\pi$}{π}}

In this section, we use decorrellation estimates and the quantity $\beta(r,R)$
to compare the crossing probabilities of quads with pairwise close boundaries,
in preparation for the scale induction in the next section. Our first estimate
is a direct use of $\beta$:

\begin{lemma}[Thickening]
  \label{lem:thickeningv2}
  Let $Q$ and $Q'$ be two quads such that $Q' \subseteq Q$, let $r>0$, and
  assume that for every $\ast \in \{\ell,r,t,b\}$, we have the proximity
  condition \[\forall x \in \partial_\ast(Q'), \quad \exists y \in
  \partial_\ast(Q), \quad d_Q(x,y) \leq r.\] Then for any boundary condition
  $\omega$ outside $Q$, \[P_{Q,\omega}[\LR(Q) \sdiff \LR(Q')] \leq 16 \beta(r,
  W(Q)/2 - r).\]
\end{lemma}

\begin{proof}
  Duality implies that the event $\LR(Q) \sdiff \LR(Q')$ is the union of two
  (nondisjoint) cases: either there is an open crossing of $Q$ and a dual
  crossing of $(Q')^\ast$, or there is an open crossing of $Q'$ and a dual
  crossing of $Q^\ast$. These two are symmetric of each other, so we will give
  an upper bound for the second one, which we will refer to as $\mathcal A$
  below, and the lemma will then follow by a union bound.

  \begin{figure}[ht]
    \centerline{\includegraphics{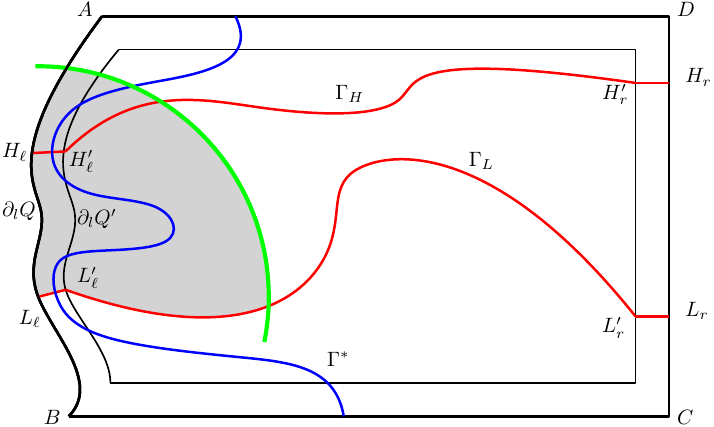}}
    \caption{The setup of the proof of Lemma~\ref{lem:thickeningv2}. The quad $Q'_L$ is shaded in grey.}
    \label{fig:thickeningv2}
  \end{figure}

  We first introduce some notation, see Figure~\ref{fig:thickeningv2}. On the
  event $\mathcal A$, denote by $\Gamma_L$ (resp.\ $\Gamma_H$) the lowest
  (resp.\ highest) open crossing of $Q'$; let $L'_\ell$ (resp.\ $L'_r$) be the
  vertex where $\Gamma_L$ meets $\partial_\ell(Q')$ (resp.\ $\partial_r(Q')$),
  $L_\ell$ (resp.\ $L_r$) a vertex on $\partial_\ell(Q)$ (resp.\
  $\partial_r(Q)$) within $d_Q$ distance at most $r$ of $L'_\ell$ (resp.\
  $L'_r$), and denote by $\Delta_{L,\ell}$ (resp.\ $\Delta_{L,r}$) a path in $Q$
  joining $L_\ell$ to $L'_\ell$ (resp.\ $L_r$ to $L'_r$) and having length at
  most $r$, which exists by assumption. The concatenation $\bar \Gamma_L$ of
  $\Delta_{L,\ell}$, $\Gamma_L$ and $\Delta_{r,\ell}$ is a path crossing $Q$
  horizontally. Introduce similar vertices and curves around $\Gamma_H$.

  Still on event $\mathcal A$, there is a vertical dual crossing $\Gamma^\ast$
  of $Q$ from $\partial_b(Q)$ to $\partial_t(Q)$. It can intersect neither
  $\Gamma_L$ nor $\Gamma_H$ (because they are open), so it has to cross
  $\Delta_{L,\ell}$ or $\Delta_{L,r}$ (possibly both), and it also has to cross
  $\Delta_{H,\ell}$ or $\Delta_{H,r}$. We focus on the case $A_{\ell,\ell}$
  where the last intersection of $\Gamma^\ast$ with $\bar \Gamma_L$ (resp.\ its
  first intersection with $\bar \Gamma_H$) belongs to $\Delta_{L,\ell}$ (resp.\
  $\Delta_{H,\ell}$), as in the figure; the other $3$ cases are handled
  similarly, and again one can conclude the proof using a union bound.

  Denote by $\Gamma^\ast_b$ (resp.\ $\Gamma^\ast_t$) the portion of $\Gamma$
  between its initial vertex on $\partial_b(Q)$ and its first hitting point on
  $\Delta_{H,\ell}$ (resp.\ between its last hitting point of $\Delta_{L,\ell}$
  and its final vertex on $\partial_t(Q)$). By definition, the diameter of
  $\Gamma^\ast$ is at least equal to $W(Q)$; moreover,
  % $\Gamma^\ast_b$ and
  % $\Gamma^\ast_t$ cover $\Gamma^\ast$ \rk{non}, up to a possible section of
  % length at most $2r$ which can only occur when $L'_\ell$ and $H'_\ell$ are
  % close; so
  at least one of the two has diameter bounded below by $W(Q)/2 - r$.
  Assume that it is $\Gamma^\ast_t$, again the other case is exactly similar.

  Now, let $D$ be the ball centered on $L_\ell$ with radius $W(Q)/2 - r$; let
  $Q_L$ be the quad lying above $\bar \Gamma_L$ in $Q$, with left boundary
  $\Delta_{L,\ell}$, bottom boundary $\Gamma_L$ and top boundary consisting in
  the arc between $L_\ell$ and $A$ along $\partial_\ell(Q)$. Note that
  $\Gamma^\ast$ is a dual crossing of $Q_L$. The intersection of $Q_L$ with $D$
  contains a quad $Q'_L$ in $\mathcal Q(r, W(Q)/2 - r)$ crossed by an initial
  section of $\Gamma^\ast_\ell$, whose left boundary is $\Delta_{L,\ell}$ and
  whose right boundary is the first arc of $\partial D \cap Q_L$ separating
  $\Delta_{L,\ell}$ from $\partial_t(Q)$ within $Q'_L$, which exists by our
  assumptions. The quad $Q'_L$ can be discovered by exploring the configuration
  outside it, so the probability that it is crossed is bounded above by
  $\beta(r, W(Q)/2 - r)$.

  Putting all the pieces together, we have $2$ choices for the crossing
  directions, $4$ choices for the topological choice of $\Gamma^\ast$, and $2$
  choices for the choice of the part of $\Gamma^\ast$ with large diameter, thus
  concluding the proof.
\end{proof}

\begin{lemma}[Separation]\label{lem:separation} If $Q$ is a quad and $Q'
  \subseteq Q$ is such that $d(Q', Q^c) \geq \ell$, then for any boundary
  condition $\omega$ on $Q^c$, \[ \big| P_{Q,\omega}[\LR(Q')] - \pi(Q') \big|
  \leq \alpha(\ell, D(Q)). \]
\end{lemma}

\begin{proof}
  We have $|P_{Q,\omega}[\LR(Q')] - \pi(Q')| \leq \tilde\alpha(Q', Q)$ so the
  upper bound follows directly from the definition of $\alpha$.
\end{proof}

\begin{lemma}[Approximation]
  For any two boundary conditions $\omega$ and $\omega'$ outside quad $Q$, and
  any $\ell < W(Q)/18$, \[ \big| P_\omega [\LR(Q)] - P_{\omega'}[\LR(Q)] \big|
  \leq 24 \beta(2\ell, W(Q)/9) + 2 \alpha(\ell,D(Q)). \]
\end{lemma}

\begin{proof}
  Choose $Q' \subseteq Q$ by quad shaving, applying Lemma~\ref{lem:shaving}.
  Then:
  \begin{align*}
    \big| P_\omega [\LR(Q)] - P_{\omega'}[\LR(Q)] \big| &\le \big| P_\omega
    [\LR(Q)] - P_\omega[\LR(Q')] \big| \\&+ \big| P_\omega [\LR(Q')] -
    P_{\omega'}[\LR(Q')] \big| \\&+ \big| P_{\omega'} [\LR(Q')] -
    P_{\omega'}[\LR(Q)] \big|.
  \end{align*}
  The first and third summand are controlled by Lemma~\ref{lem:thickeningv2},
  the second one by Lemma~\ref{lem:separation}, and the assumption on $\ell$ ensures that $W(Q)/3 - 4\ell \geq W(Q)/9$.
\end{proof}

\begin{corollary}[Oscillation]
  For any quad $Q$ and any $\ell \in \intint{2}{W(Q)/18}$, \[ \delta(Q) \leq 24
  \beta(2\ell, W(Q)/9) + 2 \alpha(\ell,D(Q)).\]
\end{corollary}

\begin{proof}
  This is a direct consequence of the previous lemma, choosing $\omega$ so that
  $P_\omega[\LR(Q)] = \pi_+(Q)$ and $\omega'$ so that $P_{\omega'}[\LR(Q)] =
  \pi_-(Q)$.
\end{proof}

The main conclusion of this section, which will be used further on to implement
our scale induction construction, can be summarized as follows:

\begin{proposition}[Optimization]
  \label{lem:epsfrombeta}
  % For any $R>0$ and any $\ell<R/24$, we have \[\epsilon(R) \leq 24 \beta(2\ell,
  % R/6 - 4\ell) + 2 \alpha(\ell,4R).\]
  For any $0<w\leq d$ and any $\ell \in \llbracket 2, w/18 \rrbracket$,
  \[\epsilon(w,d) \leq 24 \beta(2\ell, w/9) + 2 \alpha(\ell,d).\]
\end{proposition}

\begin{proof}
  % For $Q$ satisfying $W(Q) \geq R/2$ and $D(Q) \leq 4R$, and for $\ell \leq
  % R/24$, we have $\ell \leq W(Q)/12$ so we can apply the previous lemma to
  % obtain \[ \delta(Q) \leq 24 \beta(2\ell, W(Q)/3 - 4\ell) + 2 \alpha(\ell,D(Q))
  % \leq 24 \beta(2\ell, R/6 - 4\ell) + 2 \alpha(\ell, 4R). \qedhere\]

  Let $Q$ be a quad with $W(Q) \geq w$ and $D(Q) \leq d$: applying the previous
  corollary, we get \[\delta(Q) \leq 24 \beta(2\ell, W(Q)/9) + 2
  \alpha(\ell,D(Q)) \leq 24 \beta(2\ell, w/9) + 2 \alpha(\ell,d).\] Optimizing
  over all such quads leads to the claim.
\end{proof}

\section{Scale induction}
\label{sec:scale}

We are now ready to implement the scale induction constituting the core of our
argument.

\begin{lemma}[Last scale]\label{lem:lastscale} There is a universal constant
  $C>0$ such that the following happens: for every $0<\ell,\rho<R$ with $36\ell < \rho$ and every
  quad $Q$ crossing the annulus $\mathcal A(R,2R)$, \[ \pi_+(Q) \leq \max(1 -
  2^{-140} + 5 \alpha(\ell,4R) + 50 \beta(2\ell, \rho/18), \;
  \beta(\rho/2, R/2) ). \]
\end{lemma}

\begin{proof}
  Assume first that $W(Q) \leq \rho$: then there exists a cross-cut of diameter
  at most $\rho$ across $Q$, delimiting at least one quad $Q'$ crossing an
  annulus between radii $\rho/2$ and $R/2$ that any crossing of $Q$ has to
  cross. In particular, we get \[\pi_+(Q) \leq \beta(\rho/2, R/2)\] and the
  stated bound holds. From now on, we will only consider the case when $W(Q) >
  \rho$. The following argument is similar in spirit to the proof of
  Lemma~\ref{lem:thickeningv2}.

  The intersection of $Q$ with each of the circles of radii $R$ and $2R$
  consists in finitely many intervals. We call such an interval \emph{essential}
  if it is a cross-cut of $Q$ (or equivalently, if any crossing of $Q$ has to
  intersect it) and \emph{non-essential} in the opposite case. Essential
  intervals are topologically ordered and because $Q$ crosses the annulus, there
  is at least one essential interval at radius $R$ and one at radius $2R$, so
  one can always choose a pair of consecutive essential intervals, one of which
  being at each of the two radii. The portion of $Q$ between these two intervals
  is itself a quad $Q'$ crossing the annulus $\mathcal A(R,2R)$, with (up to
  symmetry) left boundary contained in the circle of radius $R$ and right
  boundary contained in the circle of radius $2R$. By Lemma~\ref{lem:ordering}, we know
  that $\pi_+(Q) \leq \pi_+(Q')$, so it is enough to prove the stated upper
  bound for $Q'$.

  By our choice of essential intervals above, there exists a crossing $\gamma$
  of $Q'$ that is contained in $\mathcal A(R,2R)$. Consider now the connected
  component $Q''$ of $Q' \cap \mathcal A(R-\rho,2R+\rho)$ which contains
  $\gamma$. It is a quad, with left and right boundaries identical to those of
  $Q'$, and top and bottom boundaries formed of portions of those of $Q'$
  concatenated with intervals of the boundaries of $\mathcal A(R-\rho,2R+\rho)$.
  It has diameter at most $2R+2\rho$ and satisfies $W(Q'') > \rho$.

  Choose $\ell>0$ and construct $Q'''$ from $Q''$ by quad-shaving. Outside an
  event of size at most $\alpha(\ell,2R+2\rho)$, the configuration $\omega$ in
  $Q$ can be coupled with a whole-plane configuration $\tilde \omega$ in such a
  way that they coincide on $Q'''$. From the annulus estimates, with probability
  at least $2^{-140}-2\varepsilon(R/2,4R)$, $\tilde \omega$ contains a closed
  circuit $\Gamma$ contained in the annulus $\mathcal A(R+\ell, 2R - \ell)$.
  Now, $\Gamma$ separates the circles of radii $R$ and $2R$, so its intersection
  with $Q'$ has to contain a cross-cut $\Gamma_0$ of $Q'$. Vertices of
  $\Gamma_0$ at distance more than $\ell$ from the boundary are closed by
  construction, so any open crossing of $Q'$ has to avoid them and hence has to
  cross some quad crossing an annulus between radii $\ell$ and $\rho$. Overall,
  this gives an upper bound on $\pi_+(Q')$ of order \[1 - 2^{-140} +
  \alpha(\ell,4R) + \beta(\ell, \rho) + 2\varepsilon(R/2,4R).\] Combining this
  with the estimate from Proposition~\ref{lem:epsfrombeta} concludes the proof.
\end{proof}

\begin{corollary}[New scale]
  \label{lem:newscale}
  For every $0<\ell,\rho,r<R$ satisfying $36\ell < \rho$, \[ \frac {\beta(r,
  2R)} {\beta(r, R)} \leq \max(1 - 2^{-140} + 5 \alpha(\ell,4R) + 50
  \beta(2\ell, \rho/18), \; \beta(\rho/2, R/2) ) \]
\end{corollary}

\begin{proof}
  Pick $\eta>0$, let $Q \in \mathcal Q(r,2R)$ such that $\pi_+(Q) \geqslant
  \beta(r, 2R) - \eta$, and let $Q_-$ be its initial fraction up to radius $R$
  and $Q_+$ be its final fraction between radii $R$ and $2R$; let $\omega$ be a
  boundary condition outside $Q$ realizing $P_{Q,\omega} \geq \pi_+(Q) - \eta$.
  Any crossing of $Q$ contains a crossing of $Q_-$ and a crossing of $Q_+$, so
  \[ \beta(r,2R) \leqslant P_{Q,\omega}[\LR(Q)] + 2\eta \leq P_\omega[\LR(Q_-)]
  \; \pi_+(Q_+) + 2\eta. \] The first factor is bounded by $\beta(r,R)$,
  Lemma~\ref{lem:lastscale} controls the second factor, so
  \[ \beta(r,2R) \leqslant \beta(r, R) \max(1 - 2^{-140} + 5 \alpha(\ell,4R) +
  50 \beta(2\ell, \rho/18), \; \beta(\rho/2, R/2) ) + 2\eta \] and since the
  inequality holds for all $\eta>0$, the lemma follows.
\end{proof}

We are now ready to state and prove our main theorem. In all that follows, we
will let
\begin{equation}\label{eq:deflambda0}
  \lambda_0 := - \frac {\log (1-2^{-141})} {\log 2} > 0.
\end{equation}

\begin{theorem}[Good up to next scale]
  \label{thm:nextscale}
  For every $0 < \lambda \leq \lambda_0$ and every $\gamma>0$ there exists a
  scale $R_0 = R_0(\gamma, \lambda)$ such that the following occurs: every
  $\gamma$-decorrelated model which is $\lambda$-good up to some scale $R \geq
  R_0$ is $\lambda$-good up to scale $2R$.
\end{theorem}

\begin{proof}
  Let $P$ be a $\gamma$-decorrelated model, let $0<r<R$ and  $0 < \lambda  \leq
  \lambda_0$ and assume that $P$ is $\lambda$-good up to scale $R$. Applying
  Corollary~\ref{lem:newscale} with $\ell = R^a$ and $\rho = R^b$ for $0<a<b<1$,
  we get
  \begin{align}
    \frac {\beta(r, 2R)} {\beta(r, R)} &\leq \max(1 - 2^{-140} + 5
    \alpha(\ell,4R) + 50 \beta(2\ell, \rho/18), \; \beta(\rho/2, R/2) ) \\
    &\leq \max(1 - 2^{-140} + 80 \gamma^{-1} R^2 \ell^{-2-\gamma} + 50 \lambda^{-1} (36\ell/\rho)^\lambda, \; \lambda^{-1} (\rho/R)^\lambda ) \\
    &\leq \max(1 - 2^{-140} + 80 \gamma^{-1} R^{2-(2+\gamma)a} + 50 \lambda^{-1} 36^\lambda R^{\lambda(a-b)}, \; \lambda^{-1} R^{\lambda(b-1)} ) \label{eq:main}
  \end{align}
  as soon as $R^{b-a} > 36$. If in addition $a$ is close enough to $1$ to
  satisfy $(2+\gamma)a > 2$, all the exponents in $R$ are negative so that for
  $R$ large enough, the RHS is at most $1-2^{-141}$. This leads to the upper
  bound
  \[\beta(r,2R) \leq (1-2^{-141}) \beta(r,R) \leq (1-2^{-141})2^\lambda
  \lambda^{-1} (r/2R)^\lambda \leq \lambda^{-1} (r/2R)^\lambda \] and concludes
  the argument.
\end{proof}

One can choose explicit values for the exponents $a$ and $b$ in the proof above
to get a quantitative bound on $R_0(\gamma,\lambda)$: letting $b=1-u$ and
$a=1-2u$ with $u = \gamma / (4 + 2\gamma + \lambda)$, all the exponents of $R$
in~\eqref{eq:main} are equal to $- \lambda u$ and the bound simplifies, showing
that one can in fact set
\[R_0(\gamma,\lambda) = \left(\frac {2^{-200}} {\gamma^{-1} + \lambda^{-1}
}\right) ^{- \frac {4 + 2\gamma + \lambda}{\lambda \gamma}}.\]

\begin{corollary}\label{cor:bootstrap} For every $0 < \lambda \leq \lambda_0$
  and every $\gamma>0$, if a $\gamma$-decorrelated model is $\lambda$-good up to
  scale $R_0(\gamma,\lambda)$, then it is $\lambda$-good at all scales.
\end{corollary}

\begin{remark}
  If a model is $\gamma$-decorrelated and $\lambda$-good at all scales, then in
  fact by the same proof, for any $\lambda'<\lambda_0$, it will be
  $\lambda'$-good at large enough scales. This shows that all good models are
  somehow ``uniformly good at large enough scales'' and is an indication that at
  least part of the perturbative argument in the last section could extend
  outside the perturbative regime.
\end{remark}

\section{Large-scale behavior of good models}

The scale induction of the previous section shows that under suitable
assumptions, models that are well decorrelated and good at a small scale are
good at all scales. It remains to show how this implies that the large-scale
qualitative topological properties of the model then match those of critical
percolation, in the sense that the box-crossing property holds and that the
one-arm probability is of polynomial order.

We strongly believe that further robust properties of percolation-like planar
models, such as quasi-multiplicativity, the values of the ``universal'' critical
exponents, and scaling relations, are valid under the same hypotheses and can be
proved in a similar vein by replacing the usual gluing arguments by uses of our
bounds on $\beta$. Implementing this programme would not be trivial, and is
beyond the scope of this article.

\bigskip

\textbf{Until the end of this section, the model $P$ is assumed to be
$\lambda$-good and $\gamma$-decorrelated.} For such a model, our previous
estimates can be made more explicit. For instance,
Proposition~\ref{lem:epsfrombeta} shows that for any $0<w\leq d$ and any $\ell
\in \llbracket 2, w/18 \rrbracket$,
\begin{equation}
\epsilon(w,d) \leq 24 \lambda^{-1} (18\ell/w)^\lambda + 2 \gamma^{-1} d^2 \ell^{-2-\gamma}.
\end{equation}
Similarly to the proof of Theorem~\ref{thm:nextscale}, applying this at $\ell =
R^{1-\mu}$ with \(\mu  = \frac {\gamma} {2+\lambda+\gamma} > 0 \) leads to the
uniform bound:

\begin{lemma}
  \label{lem:bound_on_epsilon}
  For every model which is both $\lambda$-good and $\gamma$-decorrellated, there
  exists $C$ depending only on $\gamma$ and $\lambda$ satisfying for all $R>0$
  and all $k\geq1$ the bound \[ \varepsilon(R,kR) \leq C k^2
  R^{-\gamma\lambda/(2+\gamma+\lambda)}. \] In particular, for every $k\geq1$,
  $\varepsilon(R,kR)$ vanishes as $R\to\infty$.
\end{lemma}

\begin{theorem}[No percolation]
  A model which is both good and well-decorrelated does not percolate. More
  quantitatively, if $P$ is a model which is both $\lambda$-good and
  $\gamma$-decorrelated, then there exist two constants $c>0$ and $C<\infty$
  such that for all $R>0$, \[P[0 \longleftrightarrow \partial\Lambda_R] \leq C
  R^{-c}.\]
\end{theorem}

\begin{proof}
  It is enough to prove the bound  when $R$ is a power of $2$; consider the
  collection of disjoint annuli $\mathcal A_k := \mathcal A(2^{k-1}, 2^k)$, and
  let $E_k$ be the event that $\mathcal A_k$ contains a dual circuit around the
  origin. Let $k_0$ be such that $\varepsilon(2^{k-2}, 2^{k+1}) < 2^{-142}$ for
  all $k \geq k_0$, which exists from Lemma~\ref{lem:bound_on_epsilon}. From
  Theorem~\ref{thm:annulus}, for every $k \geq k_0$, the conditional probability
  of $E_k$ conditionally on the whole configuration outside $\mathcal A_k$ is
  bounded below by $2^{-141}$. In particular, for every $k > k_0$, \[P[E_{k_0+1}
  \cup \ldots \cup E_k] \geq 1 - (1-2^{-141})^{k-k_0}\] so that $P[0
  \longleftrightarrow \partial \mathcal A_k] \leq (1-2^{-141})^{k-k_0}$. The
  second claim of the theorem follows, with
  \[ c = - \frac {\log (1-2^{-141})} {\log 2} > 0. \] The non-percolation
  follows by letting $k$ tend to $+\infty$.
\end{proof}

Notice that on the event that the box $\llbracket 0,R\rrbracket^2$ is crossed,
which has probability $1/2$, there is at least one vertex on the left-hand side,
say $x \in \{0\} \times \llbracket 0,R \rrbracket$, which is connected to
$\partial (x + \Lambda_R)$. Using a union bound and the translation invariance
of the model, this gives the lower bound \[P[0 \longleftrightarrow
\partial\Lambda_R] \geq \frac 1 {2R}\] (and this bound actually holds for any
translation-invariant and self-dual model). Together with the upper bound just
obtained, this shows polynomial decay: the probability that $0$ is connected to
distance $R$ behaves like $R^{-\Theta(1)}$.

\bigskip

\begin{theorem}[Strong BXP]
  \label{thm:oldbootstrap} If a model $P$ is both good and well-decorrelated,
  then it satisfies the strong box-crossing property: for any quad $Q$ in
  $\mathbb R^2$ with $W(Q)>0$,
  \[ \liminf_{R \to \infty} \pi_-(RQ) > 0. \]
\end{theorem}

\begin{proof}
  Choose $\lambda>0$ and $\gamma>0$ such that the moel is $\lambda$-good and
  $\gamma$-decorrellated, and first notice that the results from
  Section~\ref{sec:rsw} together with Lemma~\ref{lem:bound_on_epsilon} show that
  \begin{equation*}
    \liminf_{R \to \infty} \pi_-(4R,R) > 0 \quad \text{and} \quad \liminf_{R \to \infty} \pi(\mathcal A(R,2R)) > 0.
  \end{equation*}
  We first extend this bound to longer rectangles. From Lemma~\ref{lem:longer}
  and the definition of $\varepsilon$ we get for all $m\geqslant 3/2$ that
  \begin{equation}
    \label{eq:nextbigm}
    \pi_-((2m-1)R, R) \geqslant \frac {\pi_-(mR,R)^2} 2 -
    2 \varepsilon(R/2, 2mR).
  \end{equation}
  Using this bound repeatedly and applying Lemma~\ref{lem:bound_on_epsilon} to
  control $\varepsilon$, we then obtain the lower bound
  \begin{equation}
    \exists C>0, \quad \forall m\geq1, \quad \exists c(m)>0, \quad \forall R>0, \quad \pi_-(mR,R) \geqslant c(m) - Cm^2 R^{-\gamma\lambda/(2+\gamma+\lambda)}.
  \end{equation}
  This implies that $\liminf \pi_-(mR,R)$ is positive for all $m>0$.

  Now, consider a general topological rectangle $Q$ in the plane, with positive
  width. By a classical compactness argument, there exists a finite collection
  $(Q_i)_{1 \leq i \leq n}$ of rectangles such that
  \begin{itemize}
    \item for all $i$, $Q_i$ intersects $Q_{i+1}$ like in
      Lemma~\ref{lem:rect_to_L};
    \item $Q_i \cap Q_j = \varnothing$ as soon as $|j-i| \geqslant 2$;
    \item if for every $i$, $\gamma_i$ is a curve crossing $Q_i$, then the union
      of the $\gamma_i$ crosses $Q$.
  \end{itemize}
  These can \emph{e.g.} be obtained by dicretizing $Q$ with a square lattice of small enough mesh, choosing a self-avoiding path on that lattice which crosses $Q$, and then taking rectangular neighborhood of the straight sub-paths of that path to form the $Q_i$; see Figure~\ref{fig:genquad}.

  \begin{figure}[ht!]
    \begin{center}
      \includegraphics{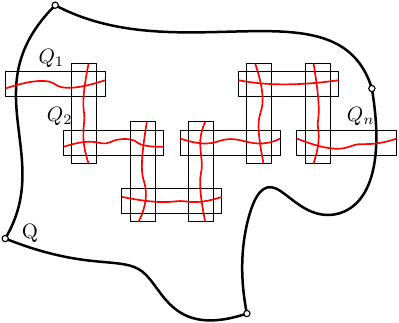}
    \end{center}
    \caption{A chain of rectangles crossing a topological rectangle. The union of the crossings $\gamma_i$ (in red) contains a crossing of $Q$.}
    \label{fig:genquad}
  \end{figure}

  For every $R>0$, the collection $(RQ_i)$ of rectangles crosses the scaled quad
  $RQ$ in the same way, and as above there exist constants $(c_i)$ and $(C_i)$
  depending only on the geometry of the $Q_i$ such that for every $i$ and every
  $R$, $\pi_-(RQ_i) \geqslant c_i - C_i R^{-\gamma\lambda/(2+\gamma+\lambda)}$.
  Performing the same gluing argument as in Lemma~\ref{lem:rect_to_L} along the
  chain of rectangles and using the previous bound on $\varepsilon$, this leads
  to a lower bound of the form \[\exists c,C>0, \quad \forall R>0, \quad
  \pi_-(RQ) \geqslant c - C.R^{-\gamma\lambda/(2+\gamma+\lambda)}\] where again
  the constants $c$ and $C$ depend only on the geometry and relative position of
  the $Q_i$ (if needed, one can check in the same way as in the beginning of
  this proof that for $n$ rectangles of aspect ratio bounded by $m$, one can
  choose $C = \mathcal O(m^2 n^2)$). This concludes the argument.
\end{proof}

\section{Interpolation}
\label{sec:interpolation}

Now, let $(P_u)_{u \in \mathbb [0,U]}$ be a one-parameter family of models, for
some $U > 0$. In addition to requiring conditions (P1)--(P3) to hold for each of
the $P_u$, assume the following:
\begin{enumerate}
\item[\newassumption{as:markov}] There exists $\rho>0$ such that, for every $u$,
  $P_u$ has the \emph{domain Markov property} with range $\rho<\infty$, in the
  sense that whenever two configurations $\omega$ and $\omega'$ coincide on a
  $\rho$-neighborhood of a finite set $D$, the probability measures
  $(P_u)_{D,\omega}$ and $(P_u)_{D,\omega'}$ are equal.
\item[\newassumption{as:continuous}] The specifications are continuous: namely,
  for every domain $D$, every configuration $\omega$ outside $D$ and every event
  $A$ measurable in the configuration inside $D$, the map $u \mapsto
  (P_u)_{D,\omega}[A]$ is continuous;
% \item[\newassumption{as:continuous2}] The map $u \mapsto P_u$ is continuous, in
%   the sense that for every event $A$ depending on finitely many vertices,
%   $P_u[A]$ depends continuously on $u$; \rk{Inutile ?}
\item[\newassumption{as:decor}] There exists $\gamma>0$ such that $P_u$ is
  $\gamma$-decorrelated for all $u \in [0,U]$;
% \item[\newassumption{as:good}] There exists $\lambda>0$ such that the model
%   $P_0$ is $\lambda$-good at all scales \rk{No this does not work}
\item[\newassumption{as:perco}\label{as:last}] The model $P_0$ is critical
  Bernoulli site-percolation.
\end{enumerate}

\begin{theorem}\label{thm:family} For every collection $(P_u)$ of models
  satisfying assumptions \ref{as:unique} to \ref{as:last}, there exist $u_0>0$
  and $\lambda>0$ such that for every $u \in [0,u_0]$, $P_u$ is $\lambda$-good
  at all scales, satisfies the box-crossing property, and does not percolate.
\end{theorem}

\begin{proof}
  For all $0<r<R$, let $\mathcal Q'(r,R)$ denote the collection of all triples
  $(Q, \partial_i, \partial_o)$ where $Q$ is a subset of $\mathcal A(r,R)$,
  $\partial_i$ a subset of $Q$ contained in $\partial \Lambda_r$ (the
  \emph{inner boundary} of $Q$) and $\partial_o$ is a subset of $Q$ contained in
  $\partial \Lambda_R$ (the \emph{outer boundary} of $Q$). The inner and outer
  boundaries of $Q$ are not required to be connected; quads contained in
  $\mathcal A(r,R)$ with $\partial_\ell(Q) \subseteq \partial \Lambda_r$ and
  $\partial_r(Q) \subseteq \partial \Lambda_R$ can be seen as special elements
  of $\mathcal Q'(r,R)$, and we will refer to generic elements of $\mathcal
  Q'(r,R)$ as \emph{generalized quads} by analogy. We make no topological
  assumption on generalized quads, but only those containing a path joining
  their inner and outer boundaries will be relevant below. Notice that the
  collection $\mathcal Q'(r,R)$ is finite (while $\mathcal Q(r,R)$ is not).

  % \bigskip

  We will use the same convention of using the same letter $Q$ to denote a
  generalized quad and its set of vertices, $\partial_i(Q)$ and $\partial_o(Q)$
  for its inner ant outer boundaries, and extend some of the notation of the
  previous sections: for every such $Q$, a crossing of $Q$ is a path in $Q$
  joining $\partial_i(Q)$ to $\partial_o(Q)$, $LR(Q)$ is the event that there
  exists an open crossing of $Q$, and for a given model, $\pi(Q)$, $\pi_+(Q)$
  and $\pi_-(Q)$ are defined correspondingly.

  In a similar way as in section~\ref{sec:scale}, any crossing of a quad $Q \in
  \mathcal Q(r,2R)$ contains a crossing of a quad $Q_- \in \mathcal Q(r,R)$ and
  a crossing of a quad $Q_+ \in \mathcal Q(R,2R)$, disjoint of each other and
  both contained in $Q$. Let $Q' \subseteq Q_+$ be the generalized quad with
  domain $Q_+ \cap \mathcal A(R,2R)$ and boundaries $Q' \cap \partial \Lambda_R$
  and $Q' \cap \partial \Lambda_{2R}$ (which is not in general a quad): any
  crossing of $Q_+$ has to contain a crossing of $Q'$, \emph{e.g.} between a
  last visit to $\partial \Lambda_R$ and a first visit to $\partial
  \Lambda_{2R}$, so $\pi_+(Q_+) \leq \pi_+(Q')$. This leads, for any model, to
  the bound
  \begin{equation}
    \label{eq:nextscalegen}
    \forall 0<r<R, \quad \frac {\beta(r,2R)} {\beta(r,R)} \leq \max_{Q' \in \mathcal Q'(R,2R)} \pi_+(Q')
  \end{equation}
  which will play the same role here as lemma~\ref{lem:newscale} above.

  % \bigskip

  Recall assumption~\ref{as:perco}: by standard Russo-Seymour-Welsh theory for
  Bernoulli percolation, the model $P_0$ satisfies the box-crossing property and
  for every $R>0$, there exists a dual circuit in the annulus $\mathcal A(R)$
  with probability bounded below by a uniform positive constant $\zeta>0$;
  section~\ref{sec:rsw} \emph{e.g.} gives $2^{-140}$ as a concrete bound, but
  the actual value is of no relevance in the context of this proof. By duality,
  this implies that the crossing probability for the annulus is bounded above by
  $1-\zeta$, and by monotonicity (which holds because $P_0$ is a product
  measure), so is the crossing probability for any generalized quad in $\mathcal
  Q'(R,2R)$.

  Now, recall $\lambda_0$ as defined by \eqref{eq:deflambda0}, choose
  \[\lambda := \min \left( \lambda_0, - \frac {\log (1-\zeta/2)} {\log 2}, 1 -
  \frac \zeta 2 \right)\] and let $R_0 = R_0(\gamma,\lambda)$ be given by
  Theorem~\ref{thm:nextscale}, with $\gamma>0$ given by
  assumption~\ref{as:decor}. There are finitely many generalized quads in the
  union of the $\mathcal Q'(R,2R)$ over $R \leq R_0$, and under $P_0$ they are
  all crossed with probability at most $1-\zeta < 1$.

  Let $Q$ be an element of $\mathcal Q'(R,2R)$; list all the (finitely many)
  configurations in a $\rho$-neighborhood of $Q$, and extend each of them
  arbitrarily to a configuration in the whole lattice, leading to a finite
  collection $(\omega_i)_{i\in I}$ of configurations. Assumption~\ref{as:markov}
  shows that for every $u$, \[\pi_{u,+}(Q) = \max_{i\in I} P_{u,Q,\omega_i}
  (LR(Q)).\] Moreover, each of the terms in this maximum depends continuously on
  $u$ by assumption~\ref{as:continuous}, therefore $\pi_{u,+}(Q)$ is itself
  continuous in $u$. By the finiteness of $\mathcal Q'(R,2R)$ there exists
  $u_0>0$ such that
  \[\forall u \in [0,u_0], \quad \forall R \leq R_0 \quad \forall Q \in \mathcal
  Q'(R,2R), \quad \pi_{u,+}(Q) \leq 1-\zeta/2.\]

  For any $0<r<R<R_0$ and any $u \in [0,u_0]$, one can iterate
  \eqref{eq:nextscalegen} over radii $2^kr$, applying it $\lfloor
  \log(R/r)/\log2 \rfloor$ times, to obtain \[\beta_u(r,R) \leq
  (1-\zeta/2)^{\lfloor \log(R/r)/\log2 \rfloor} \leq \frac
  {(r/R)^{-\log(1-\zeta/2)/\log2}} {1-\zeta/2} \leq \lambda^{-1} \left( \frac r
  R \right)^\lambda.\] In other words, all the $P_u$ for $u \in [0,u_0]$ are
  $\lambda$-good up to radius $R_0$: applying Theorem~\ref{thm:nextscale} then
  concludes the proof.
\end{proof}

\begin{remark}
  It would be very desirable to have a strengthening of Theorem~\ref{thm:family}
  where assumption~\ref{as:perco} is replaced by a weaker condition on $P_0$,
  typically stating that $P_0$ is $\lambda$-good at all scales: indeed, this
  would imply that the set of parameters $u$ where $P_u$ is nicely behaved is
  open rather than only a neighborhood of~$0$, which could open a way to a
  non-perturbative argument. Unfortunately we were not able to obtain this
  strengthening, and in fact obtaining estimates on the crossing probabilities
  of generalized quads from bounds on the $\beta(r,R)$ for a model with
  dependence on boundary conditions seems to be rather problematic.
\end{remark}

% \begin{lemma}
%   Under these assumptions, for every $r<R$, the map $u \mapsto \beta_u(r,R)$ is
%   lower semicontinuous.
% \end{lemma}

% \begin{proof}
%   Since $\beta(r,R)$ is defined as a supremum over $\mathcal Q(r,R)$, it is
%   enough to show that for every quad $Q$, the map $u \mapsto \pi_+(Q)$ is
%   continuous, so fix a quad $Q$ for the rest of the argument.

%   By the domain Markov property \ref{as:markov}, $P_\omega[\LR(Q)]$ only depends
%   on $\omega$ through its restriction to a neighborhood $Q^\rho$ of $Q$ of
%   thickness $\rho$; this neighborhood being finite, we can write \[ \pi_+(Q) =
%   \max _{\omega : Q^\rho \to \{\pm1\}} P_{u,\omega}[\LR(Q)]. \] But by the
%   continuity assumption \ref{as:continuous}, the probability under
%   $P_{u,\omega}$ of every configuration inside $Q$ depends continuously on $u$,
%   thus concluding the argument.
% \end{proof}

% \begin{proof}[Proof of Theorem~\ref{thm:family}]
%   Choose $R_0 = R_0(\lambda/2, \gamma)$ like above. Since $P_0$ is
%   $\lambda$-good, we have $\beta_0(r,R) \leq \lambda^{-1} (r/R)^\lambda$ for all
%   $r<R<R_0$, hence $\beta_0(r,R) < (\lambda/2)^{-1} (r/R)^{\lambda/2}$. Since
%   there are finitely many such pairs $(r,R)$ and all the maps $u \mapsto
%   \beta_u(r,R)$ are lower semicontinuous, we get that \[ \exists \eta>0, \forall
%   u\in(0,\eta), \forall r<R<R_0, \beta_u(r,R) < (\lambda/2)^{-1}
%   (r/R)^{\lambda/2}.\] It is now a direct consequence of
%   Corollary~\ref{cor:bootstrap} that the model $P_u$ is $\min(\lambda/2,
%   \lambda_0)$-good for all $u \in [0,\eta]$, which is what we claimed.
% \end{proof}

\begin{proposition} \label{prp:ising} Let $P_u$ denote the Ising model on the
  face-centered square lattice, with negative inverse temperature $\beta = -u$.
  Then, for $u_0$ small enough, the family $(P_u)_{u\in[0,u_0]}$ satisfies the
  assumptions \ref{as:unique}--\ref{as:last} of Theorem~\ref{thm:family}.
\end{proposition}

\begin{proof}
  We will check each assumption separately; all of them are well-known.

  \ref{as:unique} is the uniqueness of the Gibbs measure of the Ising model in
  high temperature. There are several ways to obtain it in the perturbative
  regime where $|\beta|$ is small, the most classical being Dobrushin's
  criterion~\cite{Dob68}.

  \ref{as:symm}, \ref{as:duality} and \ref{as:perco} can readily be seen from
  the definition of the specification, as does the domain Markov property
  \ref{as:markov} with range $1$ (which does hold for any Hamiltonian written as
  a sum of pairwise neighbor interactions); similarly, continuity
  \ref{as:continuous} follows directly from the fact that the interactions
  depend continuously on the parameter $\beta$.

  % \ref{as:good} is a statement about the model at $\beta=0$ which is exactly
  % critical site-percolation on the lattice. Here, boundary conditions do not
  % matter and $\beta(r,R)$ is upper-bounded by the probability that there exists
  % an open path crossing the annulus $\beta(r,R)$. This in turn is bounded by $C
  % (r/R)^c$ for universal constants $0 < c,C < \infty$ by classical
  % Russo-Seymour-Welsh results, see \emph{e.g.}~\cite{Grimmett}.

  \ref{as:decor} as we need it here is perhaps a little less standard, but by no
  means new. For instance, \emph{disagreement percolation} shows that for two
  finite domains $D' \subseteq D$, $\tilde \alpha(D',D)$ is bounded above by the
  probability, for Bernoulli percolation at parameter $p = \mathcal O(|\beta|)$,
  that there is an open path from $D'$ to $D^c$: this is exactly Corollary~1
  in~\cite{vbBM94}. Exponential decay of connectivity for subcritical
  percolation then gives a bound for this probability which is of order $|D'|
  e^{-\eta d(D', D^c)}$ for some positive $\eta$ (depending on $\beta$, but
  uniformly bounded below for $|\beta|$ small); good decorrelation follows.
\end{proof}

Theorem~\ref{thm:main} is then a direct consequence of Theorem~\ref{thm:family}
and Proposition~\ref{prp:ising}.

\bibliographystyle{siam}
\bibliography{Percolation}

\vfill

\noindent  Univ.  Grenoble  Alpes, CNRS,  Institut  Fourier,  F--38000
Grenoble, France
\end{document}